# Semi-Local Exotic Lagrangian Tori in Dimension Four

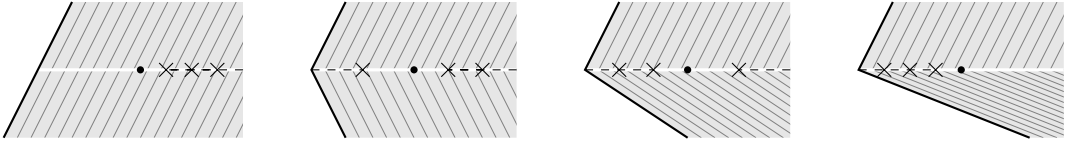

Joé Brendel, Johannes Hauber, Joel Schmitz

March 1, 2024

**Abstract**

We study exotic Lagrangian tori in dimension four. In certain Stein domains $B_{dpq}$ (which naturally appear in almost toric fibrations) we find $d + 1$ families of monotone Lagrangian tori which are mutually distinct, up to symplectomorphisms. We prove that these remain distinct under embeddings of $B_{dpq}$ into geometrically bounded symplectic four-manifolds. We show that there are infinitely many different such embeddings when $X$ is compact and (almost) toric and hence conclude that $X$ contains arbitrarily many Lagrangian tori which are distinct up to symplectomorphisms of $X$. In dimension four arbitrarily many different Lagrangian tori were previously known only in del Pezzo surfaces.

Neither the embedded tori, nor the ambient space $X$ needs to be monotone for our methods to work.

## 1 Introduction

### 1.1 Context

This paper is motivated by the symplectic classification question of Lagrangian submanifolds: for two Lagrangian submanifolds $L, L' \subset (X, \omega)$ of a symplectic manifold, can $L$ be mapped to $L'$ by a symplectomorphism of $X$? We say that $L, L'$ are **equivalent**, and write $L \cong L'$, if such a symplectomorphism exists. If it does not, we say that $L, L'$ are **inequivalent** and write $L \not\cong L'$.

An obvious necessary condition for $L$ and $L'$ to be equivalent is that they are diffeomorphic as manifolds. In this paper, we only consider Lagrangian tori. They play a special role in symplectic geometry, since they are *everywhere*. Indeed, Darboux' theorem states that every symplectic manifold locally looks like a copy of the standard symplectic vector space which contains arbitrarily small Lagrangian tori. Furthermore, Lagrangian tori naturally appear in classical mechanics, and more specifically in the study of completely integrable systems. In fact, in this paper we heavily rely on certain completely integrable systems, so-called almost toric fibrations. We use these fibrations as an auxiliary tool to study the symplectic topology of their fibres.

Recall that every Lagrangian submanifold $L \subset X$ has two characteristic classes; the **area class** $a_L : H_2(X, L) \to \mathbb{R}$ and the **Maslov class** $m_L : H_2(X, L) \to \mathbb{Z}$, which are preserved by symplectomorphisms in the sense that if there is $\varphi \in \mathrm{Symp}(X, \omega)$ with $\varphi(L) = L'$, then
$$a_L = a_{L'} \circ \Phi, \quad m_L = m_{L'} \circ \Phi, \tag{1}$$
where $\Phi = \varphi_* : H_2(X, L) \to H_2(X, L')$. Therefore, the existence of some isomorphism $\Phi$ such that (1) holds is a necessary condition for $L$ and $L'$ to be equivalent. If such a map $\Phi$ exists, we say that $L$ and $L'$ **have the same classical invariants** and write $L \sim L'$. Recall also that $L$ is called **monotone**, with monotonicity constant $C > 0$, if $a_L = C m_L$.

A full classification of Lagrangian tori up to symplectomorphism is out of reach in most ambient spaces of dimension four and above. Note that the classification is not even known for $\mathbb{R}^4$, although there is some evidence [15] that every Lagrangian torus in $\mathbb{R}^4$ is Hamiltonian isotopic to either a product torus or a Chekanov torus. Finding sets of Lagrangian tori which have the same classical invariants, but are *not* equivalent is thus a well-motivated problem. On the one hand, it sheds some light on what a classification of Lagrangian tori could look like, and on the other hand it represents an interesting occurrence of symplectic rigidity in its own right. The study of such Lagrangian tori was initiated by [11] and [17] and has attracted much interest since, see [4, 8, 13, 22, 23, 30, 34, 35], which is a non-exhaustive list. If we have a natural source of "standard Lagrangians", like product tori or toric fibres, we will say a Lagrangian is *exotic* or *symplectically knotted* if it is not equivalent to a standard Lagrangian.

## 1.2 Main results

In this paper, we study exotic Lagrangian tori in dimension four. We construct non-equivalent Lagrangian tori in certain model spaces $B_{dpq}$. If there is a symplectic embedding $B_{dpq} \supset U \to X$ of a neighbourhood $U$ of the Lagrangian skeleton of $B_{dpq}$ into a geometrically bounded symplectic four-manifold $(X, \omega)$, then we prove the images of small enough tori are non-equivalent with respect to symplectomorphisms of the new ambient space $X$.

More precisely, let $d, p, q \in \mathbb{N}$ be integers such that $d \geq 1$ and $p, q$ are coprime. For every such triple, there is a symplectic manifold $B_{dpq}$, which will serve as model space for us. We will think of $q$ as an element in $\mathbb{Z}_p/\pm 1$, since different choices of representative define symplectomorphic spaces $B_{dpq}$. For $d = 1$, $B_{1pq} = B_{pq}$ is a rational homology ball. In general, there are two complementary ways of viewing the spaces $B_{dpq}$. The first is as a Milnor fibre of the smoothing of a cyclic quotient T-singularity, see [21, Section 7.4] for more details. The second point of view, which is the one we will adopt in this paper, is that of almost toric fibrations. The spaces $B_{dpq}$ appear naturally as spaces having an almost toric base diagram with one single vertex. In other words, these spaces are the elementary building blocks of almost toric manifolds and thus very natural objects to study. Just like every vertex of a toric moment map image yields a symplectic ball embedding, every vertex of an almost toric base diagram[1] yields a symplectic embedding of a subset of $B_{dpq}$. We postpone the discussion of almost toric geometry to Sections 1.3 and 2.

**Theorem A.** *Let $(d, p, q)$ be a triple as above. Then for every $a > 0$, there are $d + 1$ monotone Lagrangian tori $T^0_{pq}(a), \ldots, T^d_{pq}(a) \subset B_{dpq}$ with monotonicity constant $\frac{a}{2}$ which have the same classical invariants, but are pairwise inequivalent, i.e.*

$$T^k_{pq}(a) \sim T^{k'}_{pq}(a), \quad T^k_{pq}(a) \not\cong T^{k'}_{pq}(a), \quad \forall k \neq k' \in \{0, \ldots, d\}. \tag{2}$$

This is proved in Section 6.

*Remark* 1.1. Some special cases of Theorem A can be deduced from [25], see Remark 6.3 for details.

**Example 1.2.** For $d = 2, p = 1, q = 0$, we have $B_{2,1,0} \cong T^*S^2$ equipped with its standard symplectic form. By Theorem A, there are three distinct types of monotone Lagrangian tori in $T^*S^2$.

1. The torus $T^0_{1,0}(a)$ is a Chekanov-type torus, i.e. it is the image of a Chekanov torus in $\mathbb{R}^4$ under a symplectic embedding $\varphi \colon \mathbb{R}^4 \supset U \to T^*S^2$.

2. The torus $T^1_{1,0}(a)$ is Clifford-type, i.e. the image of a Clifford torus under such an embedding.

3. The torus $T^2_{1,0}(a)$ is a so-called Polterovich torus, see [2].

In order to embed the exotic tori from Theorem A into a general symplectic four-manifold $(X, \omega)$, we consider certain symplectic embeddings of the form $\varphi \colon B_{dpq} \supset U \to X$, which we call **partial $B_{dpq}$-embeddings**, see Definition 2.9. Every space $B_{dpq}$ contains a configuration consisting of

1. a Lagrangian $(p, q)$-pinwheel, and

2. $d - 1$ Lagrangian spheres,

---

[1] We exclude orbifold points in our definition of almost toric fibration.

which are *visible* in the almost toric base diagram of $B_{dpq}$, see also [21, Remark 7.10]. The union of these Lagrangian subspaces is connected and forms a chain starting with the pinwheel followed by the spheres, such that every member of the chain intersects its neighbours in one point. From a symplectic point of view this union can be viewed as Lagrangian skeleton of the Stein domain $B_{dpq}$. From the algebraic geometry point of view it can be viewed as vanishing cycle of the degeneration of $B_{dpq}$ to the corresponding cyclic quotient singularity. See [19] for details on the latter point of view.

Roughly speaking, a partial $B_{dpq}$-embedding is a symplectic embedding of an open neighbourhood of the Lagrangian skeleton of $B_{dpq}$. Our second main result distinguishes the images of small enough tori $T_{pq}^k$ via partial $B_{dpq}$-embeddings into geometrically bounded symplectic four-manifolds.

**Theorem B.** *Let $\varphi : B_{dpq} \supset U \to X$ be a partial $B_{dpq}$-embedding into a geometrically bounded symplectic four-manifold $(X, \omega)$. Then there is $a_0 > 0$ such that for all $0 < a < a_0$,*

$$\varphi(T_{pq}^k(a)) \sim \varphi(T_{pq}^{k'}(a)), \quad \varphi(T_{pq}^k(a)) \not\cong \varphi(T_{pq}^{k'}(a)), \quad \forall k \neq k' \in \{0, \ldots, d\}. \quad (3)$$

*Furthermore, if $\varphi' : B_{d'p'q'} \supset U' \to (X, \omega)$ is a partial $B_{d'p'q'}$-embedding, then there is $a_0' > 0$ such that for all $0 < a' < a_0'$, if $\varphi(T_{pq}^k(a))$ can be mapped by a symplectomorphism of $X$ to $\varphi'(T_{p'q'}^{k'}(a'))$, then $a = a'$ and either $k = k' = 0$ or $(k, p, q) = (k', p', q')$.*

In other words, we distinguish not only Lagrangian tori coming from a fixed partial $B_{dpq}$-embedding, but also Lagrangian tori[2] coming from partial $B_{dpq}$-embeddings with different $(p, q)$. As a continuation of Example 1.2, the following is an immediate consequence of Theorem B and Weinstein's neighbourhood Theorem.

**Corollary 1.3.** *Let $N \subset X$ be a Lagrangian sphere in a geometrically bounded symplectic four manifold. Then every neighbourhood of $N$ contains three different types of tori, namely a Chekanov torus, a Clifford torus and a Polterovich torus, which have the same classical invariants, but are not related by symplectomorphisms of $X$.*

By [24], every Lagrangian $(p, q)$-pinwheel induces a partial $B_{1pq}$-embedding, see also [20, Definition 2.10]. Another immediate consequence of Theorem B is the following.

**Corollary 1.4.** *Let $P \subset X$ be a Lagrangian $(p, q)$-pinwheel and $P' \subset X$ a Lagrangian $(p', q')$-pinwheel in a geometrically bounded symplectic manifold $(X, \omega)$. If $(p, q) \neq (p', q')$, then there is a one-parameter family of tori in a neighbourhood of $P$ and a one-parameter family of tori in a neighbourhood of $P'$, such that its members are pairwise inequivalent up to symplectomorphisms of $X$.*

---

[2]except if $k = k' = 0$, in which case our invariant does not detect $(p, q)$

To apply Theorem B, let us now move to examples of symplectic manifolds $(X, \omega)$ which admit partial $B_{dpq}$-embeddings. For a given $X$ it is a highly non-trivial question to decide if, and for which triple $(d, p, q)$, there is a partial $B_{dpq}$-embedding. Note that the existence of a partial $B_{dpq}$-embedding in particular implies that $X$ contains a particular configuration consisting of a Lagrangian $(p, q)$-pinwheel and $d - 1$ Lagrangian spheres. In fact, for $X = \mathbb{C}P^2$, it was shown in [20], that there is a partial $B_{dpq}$-embedding if and only if $d = 1$ and $p$ is a so-called Markov number. On the other hand, we prove the following, see Proposition 7.1 for the detailed statement.

**Proposition 1.5.** *Let $(X, \omega)$ be a compact four-dimensional symplectic manifold admitting an (almost) toric structure with a base diagram having finitely many nodes. Then there is a sequence of partial $B_{d_i p_i q_i}$-embeddings into $X$ with $p_i \to \infty$.*

This is proved by considering certain sequences of almost toric mutations, see Section 7.1. Similar ideas appeared in [18]. Proposition 1.5 might be of independent interest, as it implies for example that every such manifold contains infinitely many different Lagrangian pinwheels.

Now let $X$ be as in Proposition 1.5 and let $\varphi_i : B_{d_i p_i q_i} \supset U_i \to X$ be the sequence of partial $B_{dpq}$-embeddings. For all $i$, we obtain Lagrangian tori

$$L_i^{k_i}(a_i) = \varphi_i(T_{p_i q_i}^{k_i}(a_i)), \quad \forall k_i \in \{0, \ldots, d_i\} \tag{4}$$

for all $a_i$ such that $T_{p_i q_i}^{k_i}(a_i) \subset U_i$. By the definition of partial $B_{dpq}$-embeddings, there is a segment $(0, b_i)$ on which this holds. Applying Theorem B, we obtain the following.

**Theorem C.** *Let $(X, \omega)$ be a compact four-dimensional symplectic manifold admitting an (almost) toric structure with a base diagram having finitely many nodes. For every $i \in \mathbb{N}$, there is $a_{i,0} > 0$ such that for every pair $n \neq m$, we have*

$$L_n^{k_n}(a_n) \not\cong L_m^{k_m}(a_m), \quad \forall 0 \leq k_n \leq d_n, \ 1 \leq k_m \leq d_m, \ a_n < a_{n,0}, \ a_m < a_{m,0}, \tag{5}$$

*for the tori defined in (4).*

In other words, every such $X$ contains infinitely many one-parameter families of Lagrangian tori such that no member of one family is equivalent, by a symplectomorphism of $X$, to a member of another family. The parameter $a_i$ is an area parameter in the partial $B_{dpq}$-embedding and can be viewed as the size of the corresponding torus. One drawback of this result is that the sizes of the embeddings may become arbitrarily small along the sequence, meaning $a_{i,0} \to 0$. See Remark 7.4 for a possibility how this could be improved.

As was pointed out to us by Jonny Evans, we can apply Theorem B to certain surfaces of general type studied in [18]. The key point in that paper is interpreting Mori's theory of flips in terms of almost toric base diagrams and finding a sequence of mutations, that uses only two vertices, yielding an infinite sequence of partial $B_{1p_iq_i}$-embeddings with $p_i \to \infty$. During the process, the symplectic form is deformed and thus one ends up with a non-canonical symplectic form on the surface of general type. Using Theorem B, the following is a direct corollary of [18, Theorem 1.1].

**Corollary 1.6.** *Let X be a quintic surface or a simply connected Godeaux surface. Then there is a symplectic form $\omega$ on X for which X contains infinitely many one-parameter families of Lagrangian tori the members of which are pairwise inequivalent up to symplectomorphisms of $(X, \omega)$.*

## 1.3 Methods

To construct and understand the tori considered in this paper we use almost toric fibrations, abbreviated in the following by ATFs; to distiguish those tori we use the displacement energy germ. Let us briefly discuss these methods. For a discussion and comparison to other works and methods, and in particular the use of the count of Maslov index two $J$-holomorphic disks, we refer to Section 1.4. In Section 3.5 we apply our methods to a simple example.

By the Arnold-Liouville theorem, a compact regular fibre of a completely integrable system is a Lagrangian torus. Vianna [34, 35] noticed that a certain type of completely integrable systems, called *almost toric fibrations*, is a rich source of exotic tori. We give an overview of almost toric fibrations in Section 2.1 and refer to [21] for more details. A crucial feature of ATFs is that, on a fixed symplectic manifold, there are certain operations by which one can produce new almost toric fibrations from a given one. One such operation is called *nodal slide*. A nodal slide produces a one-parameter family of almost toric fibrations on the same space, which very often produces exotic Lagrangian tori. For example, sliding a node across the base point of Clifford torus in $\mathbb{R}^4$ produces a Chekanov torus. We discuss this simple example in detail in Section 3.5.

In terms of almost toric geometry, the model spaces $B_{dpq}$ are particularly simple. Indeed, the almost toric base diagrams of $B_{dpq}$ have $d$ nodes with coinciding eigenline, see Figure 9 for two possible ATF-base diagrams that are related to each other by changing the branch cut. Subsequently sliding $d$ nodes across a given base point produces the $d + 1$ Lagrangian tori in Theorem A. The pair $(p, q)$ determines the integral affine angle between the branch cut and the boundary of the almost toric base diagram. This is discussed in detail in Section 2.3.

As mentioned before, the spaces $B_{dpq}$ are natural building blocks in the realm of almost toric fibrations. Indeed, we note the following and refer to Lemma 2.13 for a precise statement.

**Proposition 1.7.** *Every vertex of an almost toric base diagram of a symplectic manifold $(X, \omega)$ yields a partial $B_{dpq}$-embedding into X.*

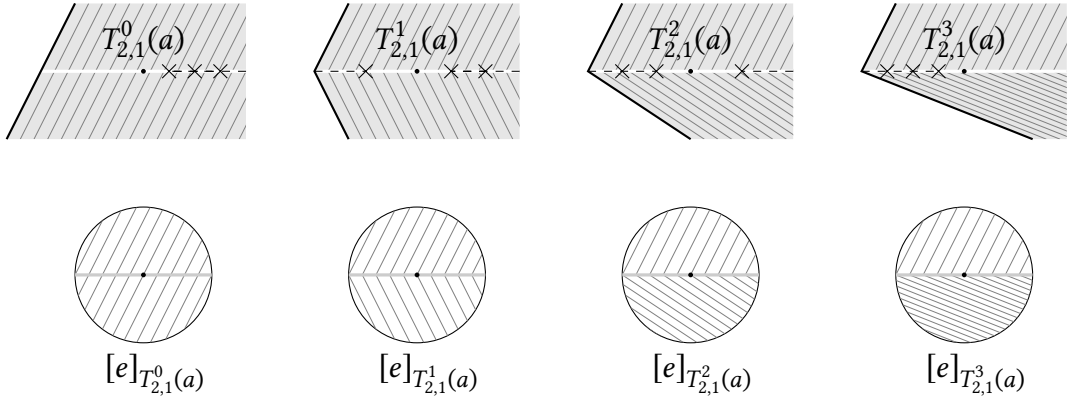

Figure 1: Upper half: Different almost toric base diagrams in which the tori $T_{2,1}^k(a)$ appear as fibres. The grey lines inside the base diagram are level sets of the displacement energy of the almost toric fibres. Lower half: Level sets of the displacement energy germs of the respective tori. For easier visibility of the level sets, the diagrams differ by a $GL(2,\mathbb{Z})$ transformation from the usual diagrams in the paper.

By iterated so-called mutations of almost toric base diagrams, one can obtain new partial $B_{dpq}$-embeddings using Proposition 1.7. This is how Proposition 1.5 is proved.

As discussed above, the tori $T_{pq}^0, \ldots, T_{pq}^d \subset B_{dpq}$ appear as fibres of different almost toric fibrations of $B_{dpq}$. More specifically, they lie on the branch cut line of its almost toric base diagram, and the number $k$ in the notation $T_{pq}^k$ of the tori is the number of nodes to one side of its base point, which means on the left by our conventions. See the upper half of Figure 1. In particular, fibres $T_{pq}^k(a)$ and $T_{pq}^{k'}(a)$ with $k \neq k'$ do not both appear as fibres in the same almost toric base diagram. Indeed, we note that the upper diagrams in Figure 1 are related by a nodal slide and a change in branch cut, for better visibility.

The invariant we use to distinguish the Lagrangian tori is the so called **displacement energy germ**. For a compact Lagrangian $L \subset (X, \omega)$ its displacement energy germ is a germ $[e]_L : H^1(L; \mathbb{R}) \to \mathbb{R} \cup \{+\infty\}$, which is derived from the displacement energy via **versal deformations**, introduced in [11]. See [7, 8, 13] for other applications of this method. Roughly speaking, the germ $[e]_L$ measures how displacement energy behaves on infinitesimal deformations of the Lagrangian $L$ in the space of Lagrangians modulo infinitesimal Hamiltonian isotopies. We carefully discuss versal deformations in Section 3.3.

If the Lagrangian in question happens to be an (almost) toric fibre, then its versal deformation is simply given by varying the base point of the fibration, see Section 3.4. Thus, in order to compute the germs $[e]_{T_{pq}^k(a)}$ (given in Theorem 6.1), we compute the displacement energy of almost toric fibres in $B_{dpq}$. This is carried out in Sections 4 and 5. See Figure 1 for a sketch of the level set of these germs. Curiously enough, we do not determine the displacement energy of the tori $T_{pq}^k(a)$ themselves.[3] Indeed, it is sufficient to know $[e]_{T_{pq}^k(a)}$ on an open dense subset in order to distinguish the tori.

In the case of the image of the tori $T_{pq}^k$ under partial $B_{dpq}$-embeddings, i.e. to prove Theorem B, the following two properties are key:

1. The displacement energy germ $[e]_L$ is well-defined, even if $L$ is not monotone. This is crucial, since, although the tori $T_{pq}^k \subset B_{dpq}$ are themselves monotone, their embedded versions into $(X, \omega)$ are not;

2. Let $\varphi \colon B_{dpq} \supset U \to X$ be a partial $B_{dpq}$-embedding. Then for small enough tori $L \subset U$, the displacement energy germ is preserved under $\varphi$,
$$[e]_{\varphi(L)} = [e]_L \circ (\varphi|_L)^*. \tag{6}$$

We say that a Lagrangian $\varphi(L)$ arising as the image of $L \subset B_{dpq}$ under a partial $B_{dpq}$-embedding has the *locality property for the displacement energy germ* if (6) holds. Note that this property is not trivial. Indeed, the germ $[e]_L$ is derived from the displacement energy of a family of tori in $B_{dpq}$, whereas $[e]_{\varphi(L)}$ is derived from the displacement energy of tori in $X$. The upper bound on the displacement energies of embedded tori, appearing as the neighbours of $\varphi(L)$, is proved by displacing them within the given partial $B_{dpq}$-embedding via the method of probes, see Section 4. The lower bound is given by the area of $J$-holomorphic curves via Chekanov's estimate Theorem 5.6, together with an area estimate on $J$-holomorphic curves leaving the partial $B_{dpq}$-embedding, see Section 5.

The locality property of the displacement energy germ was first used in [12] to distinguish embeddings of product tori and then taken up in [7] to construct exotic Lagrangian tori in Darboux charts in dimensions $\geq 6$.

## 1.4 Discussion and comparison to other works

Let us compare our results with the situation in dimensions $\geq 6$. The first-named author has recently proved that in dimension six and higher, the existence of exotic tori is generic and a purely local occurence.

**Theorem 1.8.** *[7, Theorem D] Let $(X, \omega)$ be a geometrically bounded symplectic manifold of dimension greater than or equal to six. Then every open set of $X$ contains infinitely many Lagrangian tori which have the same classical invariants, but are pairwise inequivalent up to symplectomorphisms of $(X, \omega)$.*

---

[3]Some of them are known to be non-displaceable, see Remark 4.2

The tori in Theorem 1.8 are constructed by embedding exotic tori in $\mathbb{R}^{2n}$ via Darboux charts into $X$. Although both Theorem 1.8 and our main result Theorem B are proved using the same technique, i.e. the locality property of the displacement energy germ, there are some important differences. Whereas in dimension six and above, there are known to be infinitely many inequivalent tori in $\mathbb{R}^{2n}$, which is due to [4] for monotone tori and [7] for non-monotone tori, in $\mathbb{R}^4$ it is still an open question whether there are Lagrangian tori besides product and Chekanov tori. In fact, a partial classification result in [15] suggests that there are no other Lagrangian tori in $\mathbb{R}^4$, although this is still open. This means that the approach of embedding Lagrangian tori via Darboux charts is inefficient in dimension four. Even embedding the Chekanov torus into some four-dimensional symplectic manifold $X$ by a Darboux chart does not yield a torus which is exotic in a robust sense:

*Remark* 1.9. Let $X = S^2 \times S^2$ be equipped with a non-monotone symplectic form. Then there is a Darboux embedding $\varphi \colon \mathbb{R}^4 \supset U \to X$, a Chekanov torus $L \subset U$ and a product torus $L' \subset U$ such that $\varphi(L) \cong \varphi(L')$. See also [7, Example 5.6].

This suggests that in dimension four, one should not restrict one's attention to Darboux embeddings, but look instead at embeddings of more complicated domains containing interesting tori, such as the partial $B_{dpq}$-embeddings considered here. This comes at the cost of having to construct said embeddings, for which there can be obstructions from symplectic topology, as discovered in [20], whereas Darboux embeddings exist everywhere.

Most other works on exotic Lagrangian tori have exclusively studied *monotone* exotic tori, see [4, 11, 13, 17, 23, 34, 35]. The most commonly used invariant to distinguish Lagrangian tori is a count of the Maslov index two $J$-holomorphic disks with boundary on the Lagrangian in question. Due to bubbling, this count a priori depends on the chosen $J$ in the case of non-monotone Lagrangian tori. Notable exceptions in which non-monotone exotic tori were detected are [22, 30]. In [22], the authors detect a 1-parametric family of non-displaceable tori and exoticity follows from non-displaceability in that case. In [30], the authors consider an interesting invariant $\Psi(L) \in \mathbb{R}_{>0} \cup \{+\infty\}$ based on the count of *minimal area J-holomorphic disks* as well as non-trivial input from the Fukaya category of the Lagrangian in question. In [30, Section 7] the invariant $\Psi$ is used to distinguish non-monotone Vianna-type tori in $\mathbb{C}P^2$. There is some overlap between this result and our results: From all possible mutations of the toric moment polytope of $\mathbb{C}P^2$, we obtain partial $B_{1pq}$-embeddings for all Markov numbers $p$ and thus infinitely many distinct non-monotone tori in $\mathbb{C}P^2$. The result in [30, Section 7] is quantitatively stronger,[4] since our results only apply to small tori near the ATF-vertices. On the other hand, our results apply more broadly to any partial $B_{dpq}$-embeddings into a geometrically bounded $X$.

---

[4]Although the displacement energy germ also distinguishes these tori, see [7, Example 3.7].

## 1.5 Overview

In Section 2, we discuss almost toric fibrations (Section 2.1), their fibres (Section 2.2), we introduce the model spaces $B_{dpq}$ (Section 2.3) and discuss partial $B_{dpq}$-embeddings coming from general almost toric base diagrams (Section 2.4). In Section 3, we introduce the displacement energy germ. Sections 3.1 to 3.3 give a full, self-contained account of the method of versal deformations. Versal deformations of (almost) toric fibres are discussed in Section 3.4. In Section 3.5 we explain how to combine nodal slides with computing the displacement energy germ. In Sections 4 and 5 we compute upper and lower bounds on the displacement energy of almost toric fibres of $B_{dpq}$. In Section 6, these results are combined to prove Theorems A and B. In Section 7, we prove Proposition 1.5 and discuss other examples.

## 1.6 Acknowledgements

We thank Jonny Evans for pointing out to us [18] and the example discussed in Corollary 1.6. We thank Felix Schlenk for his constant encouragement and carefully reading an earlier draft of this paper. JB acknowledges the support of the Israel Science Foundation grant 1102/20 and the ERC Starting Grant 757585.

# 2 Local model

We use almost toric geometry to define the local model spaces. In Sections 2.1 and 2.2 we give a short overview of almost toric geometry, its basic operations, which are nodal trades, nodal slides, mutations and discuss the fibres appearing in almost toric geometry. The reader familiar with the ATF-framework can skip this section. In Section 2 we introduce our local model space $B_{dpq}$, and Section 2.4 describes the relation of $B_{dpq}$ with almost toric base diagrams.

## 2.1 Almost toric fibrations and their base diagrams

Almost toric geometry is a generalization in dimension four of toric geometry. Let us briefly discuss the latter, before moving to the almost toric case. A symplectic toric four-manifold $(X, \omega, \mu)$ is a symplectic manifold $(X, \omega)$ equipped with a moment map $\mu \colon X \to \mathbb{R}^2$, generating an effective Hamiltonian $T^2$-action on $X$. The image $\mu(X) = \Delta \subset \mathbb{R}^2$ is a so called Delzant polytope and classifies $(X, \omega)$ up to $T^2$-equivariant symplectomorphisms, see [14] or [9] for more details. The moment map defines a Lagrangian torus fibration, or, equivalently,[5] a completely integrable system on $X$. The Lagrangian torus fibration defined by a toric moment map has very special properties, namely (1) its singularities, which are located over the boundary of $\Delta$, are of elliptic-regular or of elliptic-elliptic type, and (2) the moment map defines action coordinates on the interior of the moment polytope. Note that every Lagrangian torus fibration admits local action coordinates in the neighbourhood of a regular fibre, this is the classical Arnold-Liouville Theorem, but that these do not necessarily extend globally, see for example [16] or [37, 38].

---

[5]We use both terms interchangeably.

An **almost toric fibration (ATF)** is a Lagrangian torus fibration on a four-dimensional symplectic manifold whose singularities are all of elliptic-regular, elliptic-elliptic or focus-focus type. Almost toric fibrations were introduced by Symington [33], based on earlier work by Zung [36–38]. As opposed to toric manifolds, there is no globally defined Hamiltonian torus action at play, nor is there a moment map. The analogue of the moment map image is given by the so-called **almost toric base diagram**. As in the case of Delzant polytopes, every almost toric base diagram with contractible base defines a unique symplectic manifold, see [33, Corollary 5.4] and [21, Theorem 8.5].

For us,[6] the almost toric base diagram is given by a rational convex polytope $\Delta \subset \mathbb{R}^2$ in the plane decorated with nodes and branch cuts.[7] Given an almost toric fibration $F\colon X \to B \subset \mathbb{R}^2$, the corresponding base diagram is computed, roughly speaking, by attempting to find global action coordinates on the locus of regular values of $F$ in the base $B$. In the absence of focus-focus singularities, i.e. in the toric case, this succeeds and the action coordinates can be extended over the elliptic-type singularities and this procedure yields the moment map. In the presence of focus-focus singularities however, this fails. Indeed, a neighbourhood of a focus-focus singularity, whose image sits in the interior of im $F$, carries topological monodromy, meaning that the torus bundle given by neighbouring regular fibres is non-trivial. See e.g. [36] for more details. This means that action coordinates are well-defined only on the universal cover of $B\setminus\{\text{nodes}\}$. An almost toric base diagram is obtained by computing action coordinates on a fundamental domain in the universal cover, see [21, Definition 8.3] for more details. The branch cut decorations, which are represented by dashed lines, of the almost toric base diagram indicate which choice of fundamental domain was made. Equivalently, one can think of the branch cuts as curves which were removed from the regular locus of $F$ in $B$ in order to make it simply connected. On a simply connected domain, one can compute unique, meaning unique up to the usual action of $GL(2;\mathbb{Z}) \ltimes \mathbb{R}^2$, action coordinates. The image of these action coordinates yields the interior of the almost toric base diagram. The full base diagram is obtained after adding the following decorations: edges and vertices on the boundary represent elliptic-type singularities, as in the toric case, dashed lines represent the image of the branch cuts in action coordinates, and crosses represent nodes, i.e. fibres containing a focus-focus singularity. As mentioned above, a fully decorated almost toric base diagram defines a unique

---

[6]The general almost toric framework allows for much more generality, e.g. base diagrams carrying topology, non-convexity etc. but for our purposes the description we give here is sufficient.

[7]By a small abuse of notation, we denote by the symbol $\Delta$ the almost toric base diagram, as well as the subset $\Delta \subset \mathbb{R}^2$ without decorations.

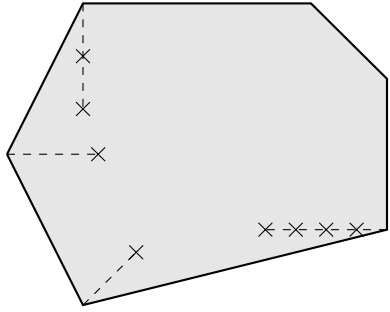

Figure 2: An almost toric base diagram. Branch cuts are drawn with dashed lines, nodes marked with crosses.

symplectic manifold.

In this paper, we will always make a special choice of almost toric base diagram. Every node in the almost toric base diagram has a unique **eigendirection**, i.e. a direction in $\mathbb{R}^2$ yielding an action coordinate which is preserved under the monodromy of the given node. Although the eigendirection is not explicitly part of the decorations of the almost toric base diagram, it can always be read off from the angle the branch cut forms at the node in question. In the current paper, we always choose the branch cut to coincide with the eigendirection. As a consequence, the almost toric base diagram "closes up", which is discussed in [21, Section 7.2], and thus yields a convex polytope. With this convention, there is an analogue of the moment map, given by a continuous map $\pi\colon X \to \Delta$, which we call **almost toric moment map**.[8] Let $\Delta_0 \subset \Delta$ be the complement of nodes and branch cuts. Assuming the branch cuts can be chosen such that $\Delta_0$ is simply connected, there are global action coordinates on its intersection with the regular locus, see [16]. Away from the branch cut and nodes, all singularities are of toric type, i.e. elliptic-regular or elliptic-elliptic, meaning that the map given by action coordinates can be smoothly extended over the singularities to yield a moment map $\pi|_{X_0}\colon X_0 \to \Delta_0$ of a Hamiltonian $T^2$-action. With our choice of branch cuts, this map has a continuous extension $\pi\colon X \to \Delta$. Let us formulate the following statement for future reference.

**Proposition 2.1.** *The restriction $\pi|_{X_0}\colon X_0 = \pi^{-1}(\Delta_0) \to \Delta_0$ is a toric moment map.*

Roughly speaking, the almost toric moment map is obtained by picking action coordinates on $X_0$, which is possible since it fibers over the simply-connected $\Delta_0$, and by continuous extension over the branch cuts.

One of the main features of almost toric geometry is that one can deform a given (almost) toric fibration $F\colon X \to \mathbb{R}^2$ to produce another one on the same symplectic manifold $X$. This freedom distinguishes the almost toric from the classical toric case. The two main operations by which one can deform an almost toric fibration are the nodal trade and the nodal slide, both of which can be handily represented by simple operations on the corresponding ATF-base diagrams. Let us briefly discuss these operations, referring to [21, Sections 8.2-8.3] for details.

---

[8]This is non-standard terminology. In contrast to "almost toric fibration", we want to empathize that it gives local action coordinates.

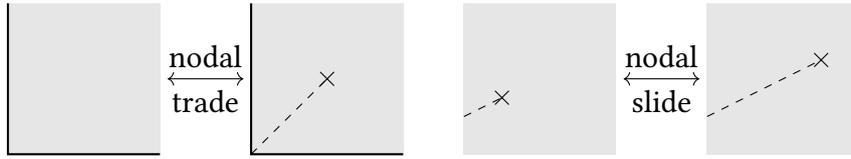

Figure 3: Nodal trade and nodal slide.

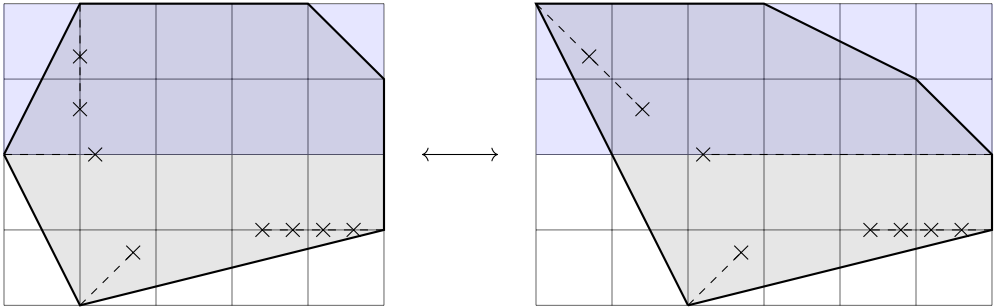

Figure 4: Changing the branch cut: Here we change the horizontal branch cut on the left, applying the corresponding shear transform (in this case $\begin{pmatrix} 1 & -1 \\ 0 & 1 \end{pmatrix}$) to the top half of the plane.

The **nodal trade** takes an elliptic-elliptic singularity and transforms it into a focus-focus singularity, while merging the two families of elliptic-regular singularities into one. In terms of the almost toric base diagram, it corresponds to replacing a toric vertex with a node and a branch cut. (See Figure 3)

The **nodal slide** consists of moving the focus-focus singularity in a non-trivial way, in terms of action coordinates, such that the node in the base diagram moves on a line in its eigendirection. (See Figure 3)

The third operation is called **changing the branch cut**, see [21, Sections 8.4] for more details. It corresponds to making a different choice of branch cut to produce the almost toric base diagram. Since, in the context of the present paper, we only consider branch cuts which coincide with the eigenline of a given node, there are only two choices of branch cuts, and the *change in branch cut*-operation is a switch from one to the other. Under such a switch, the almost toric base diagram undergoes a piece-wise linear transformation: apply the identity to one of the halves of the base diagram bounded by the eigenline and an integral shear preserving the eigenline to the other half. The branch cut is still contained in the eigenline, but switches sides. (See Figure 4)

**Proposition 2.2.** *Let $\Delta \subset \mathbb{R}^2$ be an almost toric base diagram containing a node $n \in \Delta$, lying at the end of a branch cut. Up to translation, we assume that $n = (0,0) \in \mathbb{R}^2$. Furthermore, let $v \in \mathbb{Z}^2$ be the primitive vector determining the eigenline of the node* *n* pointing away from its branch cut at n. Then the almost toric base diagram $\Delta'$ after a change in branch cut at n is obtained from $\Delta$ by applying the piece-wise linear map

$$\overline{S}_v(x) = \begin{cases} S_v(x) = x + \det(x,v)v & \text{if} \quad \det(v,x) \geq 0 \,, \\ x & \text{if} \quad \det(v,x) < 0 \,. \end{cases} \tag{7}$$

*Let $\pi$ be an almost toric moment map with image $\Delta$. Then $\pi' = \overline{S}_v \circ \pi$ is the corresponding almost toric moment map after changing the branch cut.*

The map $\overline{S}_v$ is linear on the two half-planes bounded by the eigenline of *n* and since *v* is an eigenvector of $S_v$, it is continuous. See Figure 4 for an illustration, or [21, Example 8.15] for another example.

Note that there is a difference between nodal trade and nodal slide on the one hand, and the change in branch cut operation on the other hand. The first two operations represent a change in the underlying almost toric fibration on *X*, whereas the change in branch cut does not. The change in branch cut merely corresponds to a different way of representing the almost toric fibration by a base diagram.

Together, these three operations add up to an easy-to-use visual calculus on almost toric base diagrams, all of which represent different almost toric fibrations on the same manifold. The combination of nodal slides and changes in branch cut is often called a **mutation**. For example, one can start with a toric manifold, apply nodal trades at its vertices and then perform a series of mutations. This was used in symplectic topology to great effect by Vianna [34] to exhibit infinitely many almost toric fibrations on $\mathbb{C}P^2$. The set of almost toric fibrations obtained in this way is in bijection with the set of Markov triples and each almost toric fibration has an exotic Lagrangian torus as monotone fibre. In Section 7.1 we use mutations to show that any (almost) toric manifold with (almost) toric base diagram contained in $\mathbb{R}^2$ exhibits infinitely many different almost toric fibrations.

## 2.2 Almost toric fibres

Let $F\colon X \to B \subset \mathbb{R}^2$ be an almost toric fibration with an almost toric moment map $\pi \to \Delta$. A fibre of *F* is called **almost toric fibre**. For every $x \in \Delta$, we denote by $T(x) := \pi^{-1}(x) \subset X$. If $x \in \partial\Delta$, then $T(x)$ is a point or an isotropic circle. If *x* is a node, then $T(x)$ is a Lagrangian Whitney sphere, i.e. a Lagrangian immersion of a 2-sphere with one transverse self-intersection. If *x* is neither a node, nor contained in the boundary, then $T(x)$ is a Lagrangian torus and we call such a torus a **regular ATF-fibre**.

*Remark* 2.3. A Lagrangian Whitney sphere living over the node of an ATF-base diagram can be interpreted as a limit of regular fibres, i.e. Lagrangian tori, by viewing it as a Lagrangian torus which is "pinched" along a curve on the torus. The homology class of this curve determines the topological monodromy of the torus bundle in a neighbourhood of the node.

Almost toric fibres are left unaffected by nodal slides, provided they are not contained in the segment on which the nodal slide is supported[9]:

**Lemma 2.4.** *Let $\Delta, \Delta'$ be two almost toric base diagrams which are related by nodal slides[10] that are supported in the set $\Sigma$. Let $X, X'$ be the symplectic manifolds corresponding to $\Delta, \Delta'$ and denote by $T(x) \subset X$ and $T'(x) \subset X$ their respective ATF fibres. Then for every $x \in \Delta \setminus \Sigma$, there is a symplectomorphism $\psi\colon X \to X'$ with $\psi(T(x)) = T'(x)$.*

*Proof.* This follows from the proof of [21, Theorem 8.10], by taking the set $K \subset \Delta$ therein small enough so that $x \notin K$. □

*Remark* 2.5. The claim of Lemma 2.4 fails for points which are contained in a segment on which a nodal slide is supported, even if both points correspond to regular almost toric fibres. In fact, this is one of the main sources to construct exotic Lagrangian tori. This approach was pioneered by Vianna [34, 35]. The simplest example in which this occurs is $\mathbb{R}^4 = \mathbb{C}^2$. Let $T_{Cl}(a) = S^1(a) \times S^1(a) \subset \mathbb{C} \times \mathbb{C}$ be the Clifford torus. It is the toric fibre over the base point $(a, a)$ in the standard toric structure on $\mathbb{R}^4$. To obtain an almost toric fibration with one focus-focus singularity on $\mathbb{R}^4$, perform a nodal trade at the vertex of the toric structure. By a nodal slide, move the node across $(a, a)$ to obtain the Chekanov torus $T_{Ch}(a)$ over the point $(a, a)$. See Section 3.5 for more details on the construction of the Chekanov torus, and a sketch how to distinguish it from the Clifford torus.

## 2.3 Model spaces

Let $d, p, q \in \mathbb{N}$ be integers such that $d, p \geq 1$ and $p, q$ are coprime.[11] Let $\Delta_{dpq}$ be the almost toric base diagram obtained by decorating the right half-plane

$$\Delta_{dpq} = \{(x, y) \in \mathbb{R}^2 \mid x \geq 0\}$$

with $d$ distinct nodes contained in the ray,

$$R = \{\alpha(p, q) \in \Delta_{dpq} \mid \alpha > 0\} \qquad (8)$$

and with eigendirection $(p, q)$. For a branch cut contained in $R$ and to the right of the nodes, this yields the almost toric base diagram depicted in Figure 5.

**Definition 2.6.** Let $B_{dpq}$ be the symplectic manifold having $\Delta_{dpq}$ as almost toric base diagram.

---

[9]By this we mean the line segment obtained as the union of points across which the node moves during the nodal slide.
[10]This implies that the underlying sets $\Delta, \Delta' \subset \mathbb{R}^2$ coincide.
[11]Here, we consider the pair $(1, 0)$ as coprime.

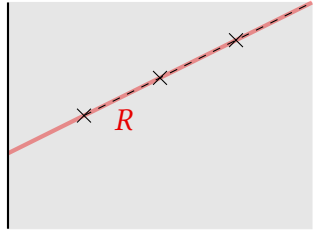

Figure 5: The almost toric base diagram $\Delta_{dpq}$ for $B_{dpq}$.

We think of $q$ as an element of $\mathbb{Z}_p/\pm 1$, as different choices of representative for $q$ are related by a $\text{GL}(2, \mathbb{Z})$ transformation, giving a symplectomorphism of the corresponding symplectic manifolds in Definition 2.6.

Such a symplectic manifolds exists, see [21, Section 7.4] and is thus unique by [21, Theorem 8.5].

*Remark* 2.7. The symplectic manifold $(B_{dpq}, \omega)$ is exact. Indeed, [21, Lemma 7.11] shows that there is a set of generators of $H_2(B_{dpq})$ represented by Lagrangian spheres, meaning that $[\omega] = 0 \in H^2(B_{dpq}; \mathbb{R})$. The space $B_{dpq}$ can be realized as the Milnor fibre of a smoothing of certain cyclic quotient singularities. For more details on this, see [21, Section 7.4] or [19].

Note that $\Delta_{dpq}$ depends on $d$ parameters $n_1 > \ldots > n_d > 0$ which correspond to the $x$-coordinate of the position of the nodes. The symplectic manifold in Definition 2.6 is independent of the position of the nodes, which is why we suppress the dependency on the $n_i$ in the notation. Indeed, any configuration of $d$ nodes on a common eigenline is related to any other such configuration by nodal slide. By [21, Theorem 8.10], this means that the space $B_{dpq}$ is well-defined up to symplectomorphism.

We turn to the almost toric fibres of $B_{dpq}$. For every $x \notin R$, we denote the corresponding almost toric fibre by $T(x) \subset B_{dpq}$. This definition is independent of the position of the nodes in $\Delta_{dpq}$ by Lemma 2.4. In case $x \in R$, the definition depends on how many nodes are to the right and to the left of $x$:

For every $a > 0$ and $k \in \{0, \ldots, d\}$, let $x_a = \left(a, \frac{aq}{p}\right) \in R$ and $\Delta_{dpq}$ be an almost toric base diagram as above for which $k$ of the $d$ nodes are to the left of the point $x_a$ and that $x_a$ is not a node (As in Figure 6).

**Definition 2.8.** Let $a > 0$ and $k \in \{0, \ldots, d\}$. We denote by $T_{pq}^k(a) = T(x_a)$ the regular almost toric fibre of a fibration whose almost toric base diagram has $k$ of the $d$ nodes to the left of $x_a$.

Again, by Lemma 2.4, this definition does not depend on the precise position of the nodes, only on the natural number $k$, i.e. the number of nodes to the left of $x_a$. In Section 6, we show that there is no symplectomorphism of $B_{dpq}$ mapping $T_{pq}^k(a)$ to $T_{pq}^{k'}(a)$ if $k \neq k'$.

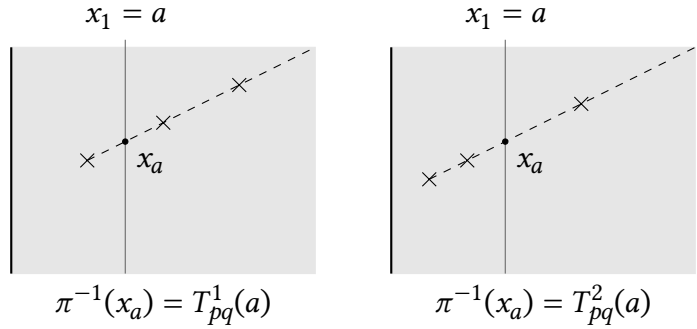

Figure 6: Modifying $\Delta_{dpq}$ by a nodal slide for Definition 2.8. The type of the Lagrangian torus changes when moving the node across its base point.

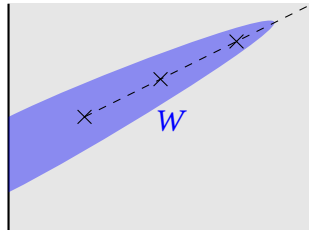

Figure 7: $W$ in the definition of a partial $B_{dpq}$ embedding.

**Definition 2.9.** A **partial $B_{dpq}$-embedding** into a symplectic manifold $(X, \omega)$ is a symplectic embedding $\varphi \colon B_{dpq} \supset U \to (X, \omega)$, where $U$

1. is connected,

2. contains a set of the form $\pi^{-1}(W)$, where $W \subset \Delta_{dpq}$ is a neighbourhood of a line segment containing the origin $(0, 0) \in \Delta_{dpq}$ and all nodes. See Figure 7.

*Remark* 2.10. In the introduction we used a slightly different definition of a partial $B_{dpq}$-embedding, namely defining a partial $B_{dpq}$-embedding as an embedding of a neighbourhood of the Lagrangian skeleton $S$ of $B_{dpq}$, as described in [21, Remark 7.10].

Indeed, the two are equivalent: Obviously Definition 2.9 gives a neighbourhood of the Lagrangian skeleton. For the other direction, we need to construct an almost toric moment map satisfying 2. in Definition 2.9. Let $L$ be the lens space $L(dp^2, dpq - 1)$, and $CL$ the cone over it. $B_{dpq} \setminus S$ is symplectomorphic to $CL \setminus \{\star\}$, where $\star$ is the orbifold point in $CL$. $CL$ has a toric moment map which is smooth except at $\star$, which can be pulled back to $B_{dpq}$ and perturbed slightly to give an almost toric moment map with the desired properties. See [19] for some more details.

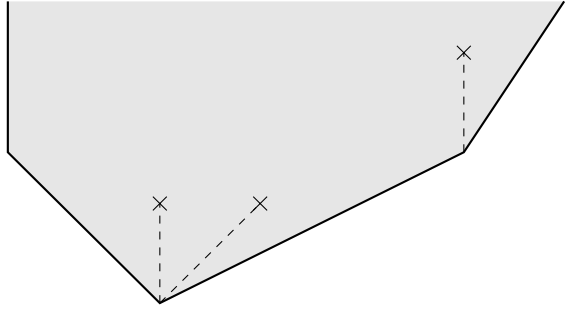

Figure 8: In our terminology, this almost toric base diagram has three vertices, one Delzant-vertex, one isolated ATF-corner and two non-isolated ATF-corners (meeting at the lower vertex).

## 2.4 Almost toric base diagrams revisited

Let $\pi : X \to \Delta$ be an almost toric moment map. In Section 2.1 we required that the branch cut of a node in $\Delta$ coincides with its eigendirection. This allows us to give a more detailed description of $\Delta$.

$\Delta$ is a rational not necessarily compact polygon[12] with two types of vertices: Either a Delzant-vertex $x \in \partial\Delta$, i.e. the primitive vectors $u_1, u_2$ pointing along the two edges emerging at the vertex form a $\mathbb{Z}$-basis of $\mathbb{Z}^2$, with no branch cut intersecting $\partial\Delta$ at $x$. See e.g. Figure 8. Or $x$ is the intersection of one or multiple branch cuts with $\partial\Delta$:

**Definition 2.11.** An **ATF-corner** of $\Delta$ is a pair $(x, dv)$ consisting of a vertex $x \in \partial\Delta$ and a branch cut intersecting $x$. By $v \in \mathbb{Z}^2$, we denote the primitive vector along the chosen branch cut which points from $x$ into the interior of $\Delta$. By $d \in \mathbb{N}$, we denote the number of nodes on that branch cut. An ATF-corner $(x, dv)$ is **isolated** if it is the only ATF-corner at the vertex $x$.

Let $(x, dv)$ be an ATF-corner of $\Delta$ and $u_1, u_2$ be the primitive vectors pointing along the edges emerging from $x$. Using Proposition 2.2, we may perform a change of branch cut at each node belonging to $(x, dv)$, thus obtaining the almost toric base diagram $\Delta' = \bar{S}_{dv}\Delta$. We may assume that $u_1$ lies in the fixed half plane of $\bar{S}_{dv}$. If $(x, dv)$ is isolated, then in $\Delta'$ there is no branch cut intersection with $\partial\Delta'$ at $x$, so either $u_1 = -S_{dv}u_2$ or $u_1$ and $S_{dv}u_2$ form a Delzant vertex. By a nodal trade we can assume that the last case does not appear. The case for non-isolated ATF-corners at $x \in \partial\Delta$ is similar, if we change the branch cut at all ATF-corners at $x$, we should map $u_1, u_2$ onto parallel vectors, such that the vertex $x \in \partial\Delta$ is mapped to the interior of an edge. Note that we can always **isolate** an ATF-corner at $x$ by changing the branch cut for all other ATF-corners at $x$.

Performing a change of branch cut for all nodes in $\Delta_{dpq}$ we get a different almost toric base diagram $\Delta'_{dpq}$ for $B_{dpq}$ with an ATF-corner $(0, d(p, q))$, see Figure 9. By

$$\Delta_{dpq} \qquad \Delta'_{dpq}$$

---
[12] An intersection of half-spaces of the form $H_{v,c} = \{x \in \mathbb{R}^2 | \langle x, v \rangle + c > 0\}$ with $v \in \mathbb{Z}^2, c \in \mathbb{R}$.

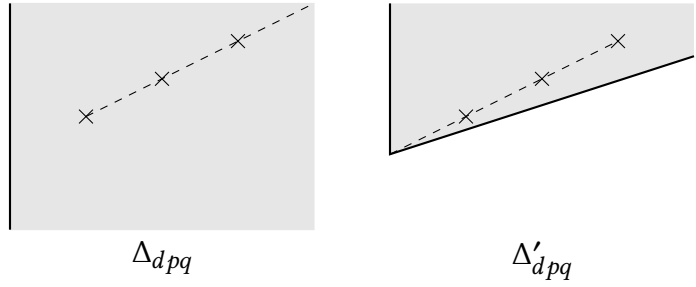

Figure 9: Changing the branch cut in $\Delta_{dpq}$ to get an ATF-corner.

a change of basis, any isolated ATF-corner $(x, dv)$ can be locally mapped by an element of $\mathrm{GL}(2, \mathbb{Z})$ to $\Delta'_{dpq}$ for an appropriate choice of $p, q$.

**Definition 2.12.** We say the ATF-corner $(x, dv)$ is of **type** $(d, p, q)$, if (perhaps after isolating it) it locally maps via a $\mathrm{GL}(2, \mathbb{Z})$ transformation into $\Delta'_{dpq}$. More precisely, there is an integral affine bijection $A : U \to V$ where $U \subset \Delta$ is a neighbourhood of the branch cut of $(x, dv)$ and $V \subset \Delta'_{dpq}$ is a neighbourhood of the branch cut of $(0, d(p, q))$.

In this sense the spaces $B_{dpq}$ are the "local models" for any almost toric fibration, as we can obtain the total space $X$ by glueing together different $B_{dpq}$ models.

**Lemma 2.13.** *If $\Delta$ has an ATF-corner of type $(d, p, q)$, then there is a partial $B_{dpq}$-embedding into $X$.*

*Remark* 2.14. We can calculate the type of an ATF-corner $(x, dv)$ as follows: Take $u_1$ to be the primitive vector pointing from $x$ along $\partial \Delta$ such that $\mathrm{int}\, \Delta$ lies to the left of $u_1$, and $u_2$ such that $u_1, u_2$ is a positively oriented $\mathbb{Z}$-basis of $\mathbb{Z}^2$. Then $p = \det(u_1, v)$ and $q = \det(u_2, v)$.

## 3 Versal deformations & nearby Lagrangians

In this section we give a self contained introduction to versal deformations. Versal defomations were previously used in [7, 8, 13] and introduced in [11].

### 3.1 Lagrangian flux

The goal of this subsection is to introduce the Lagrangian flux. This notion is described in [21, Section 2.4] in the context of Lagrangian torus fibrations, and the proofs in this section are similar to the exposition in [27, Section 10.2], which introduces the notion of symplectic flux.

Let $\mathscr{L}$ the space of Lagrangians in $X$ equipped with the $\mathscr{C}^1$-topology. See [28] for more details on the $\mathscr{C}^1$-topology on the space of Lagrangians. For the remainder of this section we fix a compact Lagrangian $L_0 \subset (X, \omega)$.

**Definition 3.1.** A **Lagrangian isotopy** based at $L_0$ is a map $\Lambda : [0,1] \to \mathscr{L}, t \mapsto \Lambda_t$ with $\Lambda_0 = L_0$ such that there is a smooth map $[0,1] \times L_0 \to X$ which maps $\{t\} \times L_0$ to $\Lambda_t$.

Note that we can always find a smooth isotopy $\varphi : [0,1] \times X \to X$ such that $\varphi_t(L_0) = \Lambda_t$. In this case we say that $\Lambda$ **is generated by** $\varphi$.

By $\mathscr{L}_{L_0} \subset \mathscr{L}$, we denote the set of Lagrangians which are Lagrangian isotopic to $L_0$. Denote by $\widetilde{\mathscr{L}}_{L_0}$ the universal cover of $\mathscr{L}_{L_0}$, that is the space of Lagrangian isotopies of $L_0$ up to endpoint preserving isotopies. Using Weinstein's Lagrangian neighbourhood theorem, one can check that $\mathscr{L}_{L_0}$ is locally simply connected, showing that $\widetilde{\mathscr{L}}_{L_0} \to \mathscr{L}_{L_0}$ is indeed the universal cover.

**Definition 3.2.** Let $\Lambda : ([0,1], 0) \to (\mathscr{L}_{L_0}, L_0)$ be a Lagrangian isotopy. Its **(Lagrangian) flux** is the map

$$H_1(L_0) \to \mathbb{R}, \quad \xi \mapsto \int_{C_\xi} \omega,$$

where $C_\xi$ is a cylinder swept out under $\Lambda$ by a representative of $\xi$ in $L_0$. We call $C_\xi$ a **flux cylinder of $\xi$ over** $\Lambda$. This map is well-defined, i.e. independent of the choice of $C_\xi$, and independent under endpoint preserving isotopies of $\Lambda$, see [32, Lemma 6.1].

Using $\mathrm{Hom}(H_1(L_0), \mathbb{R}) = H^1(L_0; \mathbb{R})$ this yields a well defined map

$$\mathrm{Flux} : \widetilde{\mathscr{L}}_{L_0} \to H^1(L_0; \mathbb{R}),$$

called the **flux map**.

**Lemma 3.3.** *A Lagrangian isotopy $\Lambda$ is generated by a Hamiltonian isotopy if and only if* $\mathrm{Flux}(\Lambda_s|_{s \in [0,t]}) = 0$, $\forall t \in [0,1]$.

*Proof.* It is shown in [32, Corollary 6.4] that the second condition is equivalent to $\Lambda$ being an exact Lagrangian isotopy. This in turn is equivalent to $\Lambda$ being generated by a Hamiltonian isotopy, see e.g. [29, Exercise 6.1.A]. □

Let $\mathscr{V} \subset \mathscr{L}_{L_0}$ be a simply-connected neighbourhood of $L_0$. Then we define the **local flux map on** $\mathscr{V}$ by setting

$$\mathrm{Flux}_{\mathscr{V}} = \mathrm{Flux}\,|_{\mathscr{V}_0}, \tag{9}$$

where we identify $\mathscr{V}$ with its lift $\mathscr{V}_0 \subset \widetilde{\mathscr{L}}_{L_0}$ containing the constant loop at $L_0$.

## 3.2 $T^*L_0$ and graphs of 1-forms

Let us look at the case $(X, \omega) = (T^*L_0, d\lambda)$, where $L_0 \subset T^*L_0$ denotes the zero section and $\lambda$ is the tautological 1-form. Recall from [27] that for any $\alpha \in \Omega^1(L_0)$, the graph $\Gamma_\alpha$ of $\alpha$ is Lagrangian if and only if $d\alpha = 0$.

**Lemma 3.4.** *Let $\alpha \in \Omega^1(L_0)$ be closed and $\Lambda_t$ be a Lagrangian isotopy from the zero section $L_0$ to the graph $\Gamma_\alpha$ of $\alpha$. Then $\mathrm{Flux}(\Lambda_t) = [\alpha]$.*

*Proof.* Let $\xi$ be a 1-cycle in $L_0$, $C_\xi$ the flux cylinder over $\Lambda$ and $\varphi_t : T^*L_0 \to T^*L_0$ an isotopy generating $\Lambda$. Then $\partial C_\xi = (\varphi_1|_{L_0})_* \xi - \xi$, on the level of homology.

The maps $(\varphi_1|_{L_0})_*, \alpha_* : H_1(L_0) \to H_1(\Gamma_\alpha)$ agree: The bundle projection induces an isomorphism $\pi_* : H_1(\Gamma_\alpha) \to H_1(L_0)$ and $\alpha_* = (\pi_*)^{-1}$. Furthermore, the family of maps $\pi \circ \varphi_t|_{L_0}$ induces the identity on homology for every $t$, since it does so at $t = 0$.

We conclude that $\partial C_\xi = \alpha_* \xi - \xi$ and compute

$$\langle \mathrm{Flux}(\Lambda_t), \xi \rangle = \int_{C_\xi} d\lambda = \int_{\partial C_\xi} \lambda = \int_{\alpha_* \xi} \lambda - \int_\xi \lambda = \int_\xi \alpha^* \lambda = \int_\xi \alpha = \langle \alpha, \xi \rangle .$$

□

**Lemma 3.5.** *Let $\alpha, \beta \in \Omega^1(L_0)$ be closed 1-forms. Their graphs $\Gamma_\alpha, \Gamma_\beta$ are Hamiltonian isotopic if and only if $[\alpha] = [\beta]$.*

*Proof.* If $[\alpha] = [\beta]$, the linear isotopy $\Lambda_t = \Gamma_{(1-t)\alpha + t\beta}$ is a Lagrangian isotopy with $\mathrm{Flux}(\Lambda_s|_{s \in [0,t]}) = 0$ for all $t \in [0,1]$. Therefore $\Lambda$ is generated by a Hamiltonian isotopy by Lemma 3.3.

Assuming that the graphs $\Gamma_\alpha$ and $\Gamma_\beta$ are Hamiltonian isotopic, we can extend a Lagrangian isotopy from $L_0$ to $\Gamma_\alpha$, for example given by $\Gamma_{t\alpha}$, by the Hamiltonian isotopy from $\Gamma_\alpha$ to $\Gamma_\beta$ to get a Lagrangian isotopy $\Lambda_t$ from the zero section to $\Gamma_\beta$. By Lemma 3.3 adding the Hamiltonian part does not change the flux map, so $[\alpha] = \mathrm{Flux}(\Lambda_t) = [\beta]$, by Lemma 3.4. □

Recall that $\mathscr{L}_{L_0}$ denotes the space of Lagrangians isotopic to $L_0$ (here it is the zero-section $L_0 \subset T^*L_0$). Let $\mathscr{L}_\Gamma \subset \mathscr{L}_{L_0}$ be the subspace of graphs of closed 1-forms. Recall also that, a priori, the flux map is defined on the universal cover of $\mathscr{L}_{L_0}$.

**Corollary 3.6.** *The subspace $\mathscr{L}_\Gamma$ contains a neighbourhood of $L_0 \in \mathscr{L}_{L_0}$ (by the definition of the $\mathscr{C}^1$-topology) and therefore $\mathrm{Flux}$ yields a well-defined (by Lemma 3.4) map*

$$\mathrm{Flux}_\Gamma : \mathscr{L}_\Gamma \to H^1(L_0; \mathbb{R}),$$

*which classifies elements in $\mathscr{L}_\Gamma$ up to Hamiltonian isotopy (by Lemma 3.5).*[13]

---

[13]In fact Ono [28, Proposition 2.3] allows to replace $\mathscr{L}_\Gamma$ by $\mathscr{L}_{L_0}$ in this statement in the case of the cotangent bundle. We do not need this stronger statement for our purposes.

## 3.3 Versal deformations

Let us return to the case of a general symplectic manifold $(X, \omega)$. Then $L_0$ admits an embedding of a Weinstein-neighbourhood $\varphi \colon T^*L_0 \dashrightarrow X$, where "$\dashrightarrow$" means that the map is only defined near the zero-section, which means that we can hope to prove a local result analogous to Corollary 3.6. Since the proof of Lemma 3.5 relies on a linear interpolation of 1-forms, we introduce the following notion.

**Definition 3.7.** A neighbourhood $\mathcal{U} \subset \mathcal{L}_{L_0}$ of $L_0$ is called **Weinstein convex** if there is a Weinstein chart $\varphi \colon T^*L_0 \dashrightarrow X$ of $L_0$ such that every $L \in \mathcal{U}$ can be written as $\varphi(\Gamma_\alpha)$ for some closed 1-form $\alpha \in \Omega^1(L_0)$, and $\mathcal{U}$ is convex with respect to the linear structure inherited from $\Omega^1(L_0)$ via $\varphi$. Such a pair $(\mathcal{U}, \varphi)$ is called **Weinstein convex pair**.

By the definition of the $C^1$-topology on $\mathcal{L}_{L_0}$, every neighbourhood of $L_0$ in $\mathcal{L}_{L_0}$ contains a Weinstein convex neighbourhood. Note that Weinstein convex sets are simply connected. Hence, the local flux map (9) makes sense on them.

**Proposition 3.8.** *Let $\mathcal{U}$ be Weinstein convex. The map* $\mathrm{Flux}_\mathcal{U}$ *classifies $\mathcal{U}$ up to Hamiltonian isotopies contained in $\mathcal{U}$. By this we mean that $L_1, L_2 \in \mathcal{U}$ are Hamiltonian isotopic by a Hamiltonian isotopy through $\mathcal{U}$ if and only if* $\mathrm{Flux}_\mathcal{U}(L_1) = \mathrm{Flux}_\mathcal{U}(L_2)$.

*Proof.* This is essentially Corollary 3.6. The fact that $\mathcal{U}$ is convex guarantees that the linear Hamiltonian isotopy of Lemma 3.5 between $L_1$ and $L_2$ is completely contained in $\mathcal{U}$. $\square$

A versal deformation is a continuous "local section of the flux map":

**Definition 3.9.** A continuous map $v_{L_0} \colon H^1(L_0; \mathbb{R}) \dashrightarrow \mathcal{L}_{L_0}$, is called a **versal deformation of** $L_0$ if for every simply connected neighbourhood $\mathcal{U} \subset \mathcal{L}_{L_0}$ of $L_0$ there exists an open neighbourhood $V$ of $0 \in H^1(L_0; \mathbb{R})$ with $\mathrm{Flux}_\mathcal{U} \circ v_{L_0}|_V = \mathrm{id}$.

By continuity of $v_{L_0}$, the neighbourhood $V$ can be chosen small enough for $v_{L_0}(V)$ to be contained in $\mathcal{U}$. Therefore, this definition makes sense.

*Remark 3.10.* It is sufficient that the condition in Definition 3.9 is satisfied for one simply connected neighbourhood $\mathcal{U} \subset \mathcal{L}_{L_0}$: If $\mathcal{V}$ is any other simply connected neighbourhood of $L_0$, then, since $\mathcal{L}_{L_0}$ is locally simply connected, there is a simply connected neighbourhood $\mathcal{W} \subset \mathcal{U} \cap \mathcal{V}$, and, since $v_{L_0}$ is continuous, restricting $V$ appropriately gives the result.

**Lemma 3.11.** *A versal deformation exists.*

*Proof.* Let $(\mathcal{U}, \varphi)$ be a Weinstein convex pair of $L_0$ and $\alpha_1, ..., \alpha_d \in \Omega^1(L_0)$ be a set of closed 1-forms inducing a basis of $H^1(L_0; \mathbb{R})$. Set

$$v_{L_0}(b) = \varphi(\Gamma_{b_1[\alpha_1]+...+b_d[\alpha_d]}), \quad b = \sum_{i=1}^d b_i[\alpha_i] \in H^1(L_0; \mathbb{R}),$$

where we have chosen a small enough domain for this definition to make sense.

Now $\mathcal{U}$ is simply connected, and by Lemma 3.4 we have $\text{Flux}_{\mathcal{U}} \circ v_{L_0} = \text{id}$. So by Remark 3.10, $v_{L_0}$ is a versal deformation. □

**Lemma 3.12.** *If $v_{L_0}, v'_{L_0}$ are two versal deformations, there is a neighbourhood $V \subset H^1(L_0; \mathbb{R})$ of $0$, such that $\forall \alpha \in H^1(L_0; \mathbb{R})$ the two Lagrangians $v_{L_0}(\alpha)$ and $v'_{L_0}(\alpha)$ are Hamiltonian isotopic.*

*Proof.* Let $\mathcal{U}$ be simply connected and $V$ such that $\text{Flux}_{\mathcal{U}} \circ v_{L_0}|_V = \text{Flux}_{\mathcal{U}} \circ v'_{L_0}|_V = \text{id}$. Additionally shrinking $V$ we may assume that $\mathcal{U}$ is Weinstein convex. Then for $\alpha \in V$, by Proposition 3.8 the result follows since $\text{Flux}_{\mathcal{U}}(v_{L_0}(\alpha)) = \text{Flux}_{\mathcal{U}}(v'_{L_0}(\alpha))$. □

Colloquially we might say that a versal deformation parametrizes all nearby Lagrangians up to local Hamiltonian isotopy by a neighbourhood of zero in $H^1(L_0; \mathbb{R})$.

Let $I: \mathscr{L} \to A$ be a symplectic invariant of Lagrangian submanifolds, meaning that $I(\varphi(L)) = I(L)$ for all $\varphi \in \text{Symp}(X, \omega)$. Any versal deformation $v_{L_0}$ gives rise to a function $I \circ v_{L_0}$. By Lemma 3.12 and symplectic invariance of $I$, the germ of this function is independent of the choice of versal deformation.

**Definition 3.13.** Let $L_0 \subset X$ be a compact Lagrangian, $v_{L_0}$ a versal deformation thereof and $I$ a symplectic invariant of Lagrangian submanifolds. By abuse of terminology, the germ of $I \circ v_{L_0}$ at $0 \in H^1(L_0; \mathbb{R})$ is called the **germ of $I$ at $L_0$**. It is denoted by
$$[I]_{L_0} : H^1(L_0; \mathbb{R}) \dashrightarrow A \, .$$

**Proposition 3.14.** *If $\varphi$ is a symplectomorphism, then*
$$[I]_{\varphi(L_0)} = [I]_{L_0} \circ \varphi^* \, ,$$
*where $\varphi^*$ denotes the induced map on cohomology.*

*Consequently, if there is no isomorphism $\Phi: H^1(L'; \mathbb{Z}) \to H^1(L; \mathbb{Z})$ such that $[I]_L \circ \Phi = [I]_{L'}$, there is no symplectomorphism $\varphi$ with $\varphi(L) = L'$.*

*Proof.* Let $v$ be a versal deformation of $L_0$. Then $v_\varphi = \varphi \circ v \circ \varphi^*$ is a versal deformation of $\varphi(L_0)$. We have the following commutative diagram:

$$\begin{array}{ccc} H^1(L_0; \mathbb{R}) & \xleftarrow{\varphi^*} & H^1(\varphi(L_0); \mathbb{R}) \\ \downarrow{v} & & \downarrow{v_\varphi} \\ \mathscr{L} & \xrightarrow{\varphi} & \mathscr{L} \\ & \searrow{I} \quad \swarrow{I} & \\ & A & \end{array}$$

where the top square commutes by definition of $v_\varphi$, and the bottom triangle commutes since $I$ is invariant under symplectomorphisms. □

## 3.4 Versal deformations of (almost) toric fibres

Let $\mu: U \to \mathbb{R}^n$ be the moment map of an effective Hamiltonian $T^n$-action defined on a subset $U \subset (X, \omega)$ of a symplectic manifold. Such a local moment map can be constructed, for example, in the neighbourhood $U$ of a regular fibre of a completely integrable system by the Arnold–Liouville theorem. In that case, $\mu$ corresponds to action coordinates. Let $L = \mu^{-1}(x_0)$ be a regular fibre of $\mu$. Since $L$ is diffeomorphic to a torus and $H^1(L; \mathbb{R}) \cong \mathbb{R}^n$, there exists an $n$-dimensional versal deformation. It is therefore natural to look at neighbouring fibres, i.e. Lagrangian tori of the form $\mu^{-1}(x_0 + a)$ for small $a \in \mathbb{R}^n$. By Definition 3.9 neighbouring Lagrangians in a versal deformation are parametrized by their Lagrangian flux. In fact, the local moment map coincides with the flux map, restricted to fibres $\mu^{-1}(x)$, see for example [21, Lemma 2.15]. In the language of versal deformations, this observation yields the following.

**Proposition 3.15.** *Let $L = \mu^{-1}(x_0)$ be a regular fibre of a (locally defined) toric moment map $\mu$. Then*

$$v_L : H^1(L; \mathbb{R}) \cong \mathbb{R}^n \dashrightarrow \mathscr{L}_L, \quad v_L(a) = \mu^{-1}(x_0 + a) \tag{10}$$

*is a versal deformation of $L$.*

The components $\mu_i$ of the moment map $\mu$ generate $n$ closed orbits in $L$, which yield a natural identification $H_1(L) \cong \mathbb{Z}^n \subset \mathbb{R}^n$. Dually, we obtain an identification $H^1(L; \mathbb{R}) \cong \mathbb{R}^n$ used in (10).

## 3.5 How to distinguish tori using the displacement energy germ

Since this will be the main technique in the proof of Theorems A and B, let us briefly outline how to combine almost toric fibrations with versal deformations to distinguish Lagrangian tori.

Let $L \subset (X, \omega)$ be a Lagrangian torus which we suspect to be exotic and let $I(\cdot)$ be a symplectic invariant of Lagrangian submanifolds as in the discussion surrounding Definition 3.13. All exotic tori known to the authors have the property that most neighbouring tori are themselves Hamiltonian isotopic to standard (non-exotic) tori, e.g. product tori in $\mathbb{R}^{2n}$ and toric fibres in toric manifolds. Therefore, the germ $[I]_L$ can be computed on an open dense set of the origin of $H^1(L; \mathbb{R})$ in the following two steps:

(1) compute $I(\cdot)$ on the standard tori of $X$,

(2) identify the standard torus to which a given member $v_L(a)$ of the versal deformation is Hamiltonian isotopic to.

We called this property the *instability of the exotic torus L* in [7] and refer to [7, Sections 1.6, 5.1] for a further discussion of instability. In our case, we use the displacement energy $I(\cdot) = e(\cdot)$ as an invariant to distinguish Lagrangian tori, see Section 5 for a discussion of $e(\cdot)$. Since $e(\cdot) \in \mathbb{R} \cup \{+\infty\}$, we obtain a germ $[e]_L : H^1(L;\mathbb{R}) \cong \mathbb{R}^n \dashrightarrow \mathbb{R} \cup \{+\infty\}$.

**Definition 3.16.** Let $L, L' \subset (X, \omega)$ be Lagrangian tori. We say that their displacement energy germs $[e]_L, [e]_{L'}$ are **equivalent**, written $[e]_L \cong [e]_{L'}$, if there is an isomorphism
$$\Phi : H^1(L';\mathbb{Z}) \to H^1(L;\mathbb{Z}), \tag{11}$$
such that $[e]_L \circ \Phi$ and $[e]_{L'}$ agree on an open neighbourhood of zero in $H^1(L';\mathbb{R})$.

If they agree only on an open dense set intersecting each neighbourhood of 0, we say they are **roughly equivalent**, and write $[e]_L \sim [e]_{L'}$.

Note that $[e]_L \not\sim [e]_{L'}$ implies $[e]_L \not\cong [e]_{L'}$.

In examples, we work with identifications $H^1(L';\mathbb{R}) \cong H^1(L;\mathbb{R}) \cong \mathbb{R}^n$ induced by a choice of basis $H_1(L) \cong H_1(L') \cong \mathbb{Z}^n$, in which case the isomorphism (11) corresponds to an automorphism $\Phi \in \mathrm{GL}(n;\mathbb{Z})$. The following is an immediate consequence of Proposition 3.14.

**Proposition 3.17.** *Let $L, L' \subset X$ be Lagrangian tori with $[e]_L \not\cong [e]_{L'}$. Then there is no symplectomorphism of $X$ mapping $L$ to $L'$.*

Computing the displacement energy germ is particularly easy for almost toric fibres. Firstly because every regular almost toric fibre has a natural versal deformation by Proposition 3.15 and secondly because in an almost toric fibration, it is particularly easy to identify the standard torus to which a versal deformation is Hamiltonian isotopic to. Compare step (2) in the above outline. Let us illustrate this by means of the basic example of the Chekanov torus in $\mathbb{R}^4 = \mathbb{C}^2$.

Let $\mu \colon \mathbb{C}^2 \to \Delta \subset \mathbb{R}^2$ be the standard moment map on $\mathbb{C}$ given by $\mu(z_1, z_2) = \left(\pi|z_1|^2, \pi|z_2|^2\right)$, with the moment polytope $\Delta \subset \mathbb{R}^2$ being the upper right quadrant. Its toric fibres are the so-called product tori, meaning products of standard circles in the plane. Using similar techniques as in Sections 4 and 5 we find that the displacement energy of product tori is

$$e(T(x)) = \min\{x_1, x_2\}, \quad x = (x_1, x_2) \in \mathrm{int}(\Delta). \tag{12}$$

This yields the answer to Step (1) in the above outline.

Let $x_0 > 0$. As a warm-up to the case of the Chekanov torus, let us compute the displacement energy germs of the product tori $T^1_{1,1}(x_0) = T(x_0, x_0)$ and $T = T(x_0, x_2)$ with $x_2 > x_0$. The former is the monotone Clifford torus, the notation $T^1_{1,1}$ of which is coherent with the one used in Definition 2.8, and the latter is a non-monotone product torus. We use Proposition 3.15 to find that $b \mapsto T(x + b)$ is a versal deformation of $T(x)$ for any $x \in \operatorname{int} \Delta$. Together with (12), this yields the displacement energy germs

$$[e]_{T^1_{1,1}}(b) = x_0 + \min\{b_1, b_2\}, \quad [e]_T(b) = x_0 + b_1, \tag{13}$$

where $b \in H^1(T; \mathbb{R}) \cong H^1(T^1_{1,1}; \mathbb{R}) \cong \mathbb{R}^2$ using the obvious identifications. By Proposition 3.14, this shows that $T$ and $T^1_{1,1}$ are distinct up to symplectomorphism of the ambient space,[14] although they both have displacement energy $x_0 > 0$. See Figure 10 for the level sets of $x \mapsto e(T(x))$ of the corresponding germs.

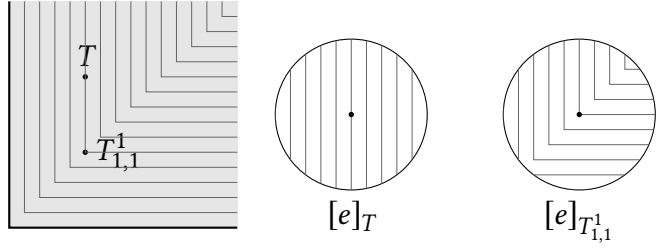

Figure 10: Distinguishing the tori $T$ and $T^1_{1,1}$

Let $\pi \colon \mathbb{C}^2 \to \Delta$ be the almost toric moment map obtained from the toric case discussed above by performing a nodal trade at the vertex of $\Delta$. One way of defining the Chekanov torus $T^0_{1,1}(x_0) \subset \mathbb{C}^2$, where $x_0 > 0$ is the area parameter, is as almost toric fibre $\pi^{-1}(x_0, x_0)$. Note that this definition works for any $x_0 > 0$, as using a nodal slide one can arrange for the position of the node to be $(y, y)$ with $y > x_0$. Let $x \in \operatorname{int} \Delta$ be a point not contained in the branch cut. Then, by Lemma 2.4, we find that $\pi^{-1}(x)$ can be mapped by a symplectomorphism to $T(x)$. Therefore $e(\pi^{-1}(x)) = \min\{x_1, x_2\}$ for such almost toric fibres. This corresponds to Step (2) in the outline above. See the upper middle diagram in Figure 11 for the level sets of $x \mapsto e(\pi^{-1}(x))$. The displacement energy of the Chekanov tori on the branch cut cannot be determined in this way, as they are not symplectomorphic to product tori. In order to apply Proposition 3.15, we use Proposition 2.2 to change the branch cut in order for the Chekanov tori not to lie on the branch cut. We find

$$[e]_{T^0_{1,1}}(b) \sim x_0 + b_1, \tag{14}$$

by which we mean equality on the open dense subset obtained as the complement of the line formed by Chekanov tori, see Definition 3.16 and Figure 11. This proves that $T^0_{1,1}$ is exotic, meaning that it cannot be mapped to a product torus by a symplectomorphism. Indeed, $T^0_{1,1}$ is monotone and hence distinct from any $T(x)$ with $x_2 > x_1$. The displacement energy germs (13) and (14) together with Proposition 3.17 show that $T^0_{1,1}$ and $T^1_{1,1}$ are distinct.

---

[14]This also follows from the fact that $T^1_{1,1}$ is monotone, whereas $T$ is not.

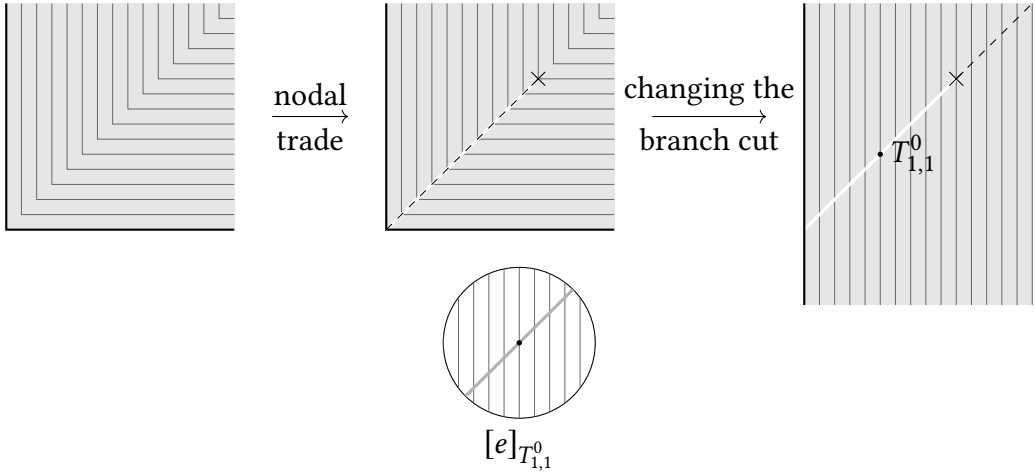

Figure 11: Constructing the Chekanov torus $T_{1,1}^0$.

For further examples of distinguishing Lagrangian tori using versal deformations see [6–8, 11–13]. For a detailed expository introduction to the Chekanov torus, see [5].

# 4 Upper bound on displacement energy: probes

Let $T(x)$ in $B_{dpq}$ be a regular almost toric fibre with $x \in \Delta_{dpq} \setminus R$, with $R$ being the eigenline of the mondromy as in (8). We use the method of probes to compute an upper bound on the displacement energy of $T(x)$. Probes were introduced by McDuff [26] to prove the displaceability of toric fibres, see also [1]. In [8], they were used to compute an upper bound on the displacement energy of toric fibres.

Our main observation is that, away from nodes and branch cuts, almost toric systems are toric, compare also Proposition 2.1, and thus the method of probes applies. For basic definitions and statements about probes, we refer to [26], especially [26, Definition 2.3, Lemma 2.4].

Recall $\Delta_{dpq}$ from Figure 5.

**Lemma 4.1.** *Let $x \in \text{int}\,\Delta_{dpq} \setminus R$. Then the displacement energy of the regular almost toric fibre $T(x)$ satisfies*
$$e(T(x)) \leq x_1,$$
*where $x_1$ is the first component of $x \in \Delta_{dpq} \subset \mathbb{R}^2$.*

*Proof.* We will use the probe $P = \{(s, x_2) \mid s \geq 0\} \subset \Delta_{dpq}$, where $x_2$ is the second component of $x \in \Delta_{dpq}$. It intersects the boundary $\partial \Delta_{dpq}$ integrally transversally in the sense of [26, Section 2.1]. However, note that $P$ may intersect nodes or the branch cut, depending on the position of nodes in $\Delta_{dpq}$. Since $x \notin R$, we can, by Lemma 2.4, carry out a nodal slide such that $P$ is disjoint from the nodes and the branch cut. By Proposition 2.1, the almost toric moment map is toric on the preimage of a neighbourhood of $P$ and thus $P$ is a probe of infinite length. See Figure 12.

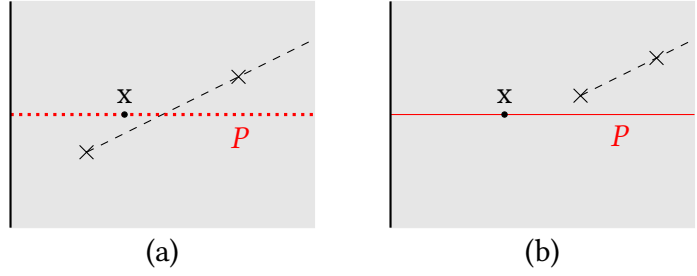

Figure 12: (a) Desired probe path $P$ obstructed by branch cut.
(b) Cleared probe path using nodal slide.

By [8, Proposition 3.4], the almost toric fibre $T(x)$ is displaceable and satisfies $e(T(x)) \leq x_1$. □

*Remark* 4.2. The above argument fails for some almost toric fibres over the ray $R$, namely for all[15] $T_{pq}^k(a)$ with $k > 0$. In fact, the tori $T_{pq}^k$ have been shown to have non-vanishing Lagrangian Floer homology, which implies that they are non-displaceable, by Lekili–Maydanskiy [25] in the following two cases (see also Remark 1.1):

1. $p = q = 1$, arbitrary $d$, and $k > 1$;

2. $d = 1$, arbitrary $p, q$, and $k = 1$.

Indeed, for the first claim, note that $B_{d,1,1}$ appears as Milnor fibre of the smoothing of the $A_{d-1}$-singularity and thus [25, Proposition 2.20] implies the claim. For the second claim, see [25, Proposition 3.6]. In the special case $B_{2,1,1} \cong T^*S^2$ of 1., the torus $T_{1,1}^2(a)$ corresponds to the so-called *Polterovich torus*, which had been shown to be non-displaceable earlier by Albers–Frauenfelder [2].

In light of this remark it seems reasonable to conjecture that all $T_{pq}^k(a) \subset B_{dpq}$ with $k > 0$ (with the exception of $p = q = k = 1$) are non-displaceable.

**Proposition 4.3.** *Let $\varphi \colon B_{dpq} \supset U \to (X, \omega)$ be a partial $B_{dpq}$-embedding. Then there exists a neighbourhood $V \subset \Delta_{dpq}$ of the origin such that for all $x \in V$, we have $e(\varphi(T(x))) \leq x_1$.*

*Proof.* We displace $T(x)$ by a Hamiltonian diffeomorphism with support in $U \subset B_{dpq}$. The proof is essentially the same as in Lemma 4.1, with the exception that the probe is not infinitely long. By Definition 2.9, the set $U$ contains the preimage $\pi^{-1}(W)$ of a neighbourhood $W$ of a line segment containing the origin and all nodes. We consider all horizontal probes contained in $W$ and let $V$ be the set of points which are displaceable by such probes. Recall from [26, Section 2.1] that this set consists of all points lying on the left half of a probe of this type. Thus $V$ is a non-empty neighbourhood of the origin. The displacement energy is estimated by the same argument as in the proof of Lemma 4.1. □

---

[15]With the exception of the case $p = q = k = 1$, which corresponds to the Chekanov torus and is displaceable by a probe.

# 5 Lower bound on displacement energy: minimal J-holomorphic curves

We want to use a classical theorem by Chekanov [10], that gives a lower bound on the displacement energy of Lagrangian submanifolds in certain symplectic manifolds. To state the theorem we recall the central notions contained in the theorem.

Let $(X, \omega)$ be a symplectic manifold. Define $\mathscr{H}(X) := \mathscr{C}_c^\infty(X \times [0,1], \mathbb{R})$, the space of compactly supported, time-dependent Hamiltonian functions on $X$. For $H \in \mathscr{H}(X)$ let

$$\|H\| := \int_0^1 \left( \sup_{x \in X} H_t(x) - \inf_{x \in X} H_t(x) \right) dt$$

be the **Hofer norm of** $H$.

**Definition 5.1.** For a subset $A \subseteq X$ of a symplectic manifold $(X, \omega)$ define the **displacement energy** of $A$ to be

$$e(A) := \inf\{\|H\| \mid H \in \mathscr{H}(X) \text{ such that } \varphi_1^H(A) \cap A = \emptyset\},$$

where $\varphi_1^H$ denotes the Hamiltonian diffeomorphism generated by $H$ and the infimum over the empty set is defined to be infinity.

Chekanov's theorem only holds for a specific class of symplectic manifolds, namely geometrically bounded ones.

**Definition 5.2.** Let $(X, \omega)$ be a symplectic manifold. We call an almost complex structure $J$ on $(X, \omega)$ **geometrically bounded**[16] if there is a Riemannian metric $g$ on $X$ such that

1. $g$ is complete,

2. there are constants $C_1, C_2 > 0$ such that $\omega(JV, V) \geq C_1 \|V\|_g^2$ and $|\omega(V, W)| \leq C_2 \|V\|_g \|W\|_g$ for all tangent vectors $V, W$, and

3. the sectional curvature of $g$ is bounded from above and the injectivity radius is bounded from below.

---

[16]Sometimes the word "tame" is used with similar meaning. See e.g. [31, Definition 4.1.1].

Denote the space of such almost complex structures by $\mathcal{J}_{\text{geo}}(X, \omega)$. The symplectic manifold $(X, \omega)$ is called **geometrically bounded** if there is a geometrically bounded $J$.

*Remark* 5.3. Closed symplectic manifolds are geometrically bounded, and the same holds for a large class of symplectic manifolds that look standard at infinity, e.g. cotangent bundles of compact symplectic manifolds and symplectizations of compact contact manifolds. For more details on this notion see also [3, Definition 2.2.1]. In particular, $B_{dpq}$ is geometrically bounded because outside a compact subset it is symplectomorphic to one end of the symplectization of a certain lens space. For more details see [19] and [21].

*Remark* 5.4. It follows from the definition that if $J_0 \in \mathcal{J}_{\text{geo}}(X, \omega)$ is a geometrically bounded almost complex structure, $K \subset M$ compact, and $J$ is a tame almost complex structure such that $J = J_0$ on $M \setminus K$, then $J$ is also geometrically bounded.

With these definitions in place we define the two quantities that play a crucial role in Chekanov's theorem.

**Definition 5.5.** Suppose that $L$ is a compact Lagrangian submanifold in a symplectic manifold $(X, \omega)$ that is geometrically bounded. Furthermore, let $J \in \mathcal{J}_{\text{geo}}(X, \omega)$ be a geometrically bounded almost complex structure. Define

$$\sigma_D(X, L, J) := \inf\left\{\int_D u^*\omega \;\middle|\; \begin{array}{c} u \colon (D, \partial D) \to (X, L) \\ \text{non-constant } J\text{-holomorphic disc} \end{array}\right\}$$

$$\sigma_S(X, J) := \inf\left\{\int_{S^2} u^*\omega \;\middle|\; \begin{array}{c} u \colon S^2 \to X \\ \text{non-constant } J\text{-holomorphic sphere} \end{array}\right\}.$$

These two quantities may be equal to infinity, if the infimum is taken over an empty set. We have that $\sigma_D, \sigma_S > 0$ in the case that $X$ is geometrically bounded. This is proven, for example, in [31, Proposition 4.3.1 (iii) and Proposition 4.7.2 (iii)].

**Theorem 5.6** (Chekanov 1998 [10]). *Let $L$ be a compact Lagrangian submanifold of a geometrically bounded symplectic manifold $(X, \omega)$. For every geometrically bounded almost complex structure $J$, the displacement energy satisfies*

$$e(L) \geq \min\{\sigma_D(X, L, J), \sigma_S(X, J)\}.$$

Let $(M, \omega, J)$ be a geometrically bounded connected symplectic manifold with geometrically bounded almost complex structure $J$, denote by $g$ a Riemannian metric as in Definition 5.2, and let $u \colon \Sigma \to M$ be a $J$-holomorphic curve, where $\Sigma$ is a compact Riemann surface.

Then Gromov's Monotonicity Lemma (e.g. [31, Proposition 4.3.1 (ii)]) gives us the following statement:

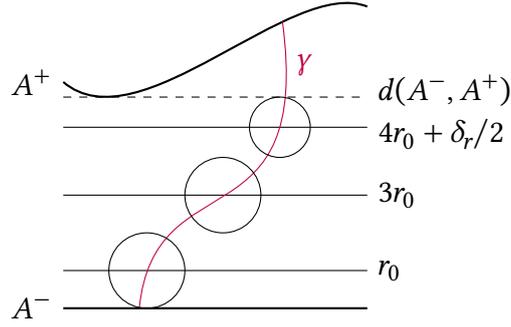

Figure 13: Placing disjoint balls on $u$ between $A^-, A^+$

**Lemma 5.7** (Monotonicity). *There exist constants[17] $C > 0, r_0 > 0$, that depend only on $(M, \omega, J)$, such that for all $x \in \mathrm{im}(u)$, $B_r(x)$ an open ball with respect to $g$ with $r < r_0$ and $u(\partial \Sigma) \subset M \setminus B_r(x)$, we have*

$$\int_u \omega \geq Cr^2 \ .$$

As a simple corollary we get that a $J$-curve crossing a certain distance, measured with respect to the Riemannian metric $g$, also must have a certain amount of area:

**Corollary 5.8.** *Suppose that $A^+, A^- \subset M$ are non-empty, that $\Sigma$ is a connected Riemann surface with $\partial \Sigma = \partial \Sigma^+ \sqcup \partial \Sigma^-$, with both $\partial \Sigma^\pm$ non-empty, and that $u : \Sigma \to M$ is a $J$-holomorphic curve with $u(\partial \Sigma^\pm) \subset A^\pm$.*

*Let $\delta = d(A^+, A^-)$, where $d$ denotes the metric induced by $g$, $N = \left\lfloor \frac{\delta}{2r_0} \right\rfloor$ and $\delta_r = \delta - 2r_0 N$. Then*

$$\int_u \omega \geq C(Nr_0^2 + \frac{\delta_r^2}{4}) \ .$$

*Proof.* Pick a path $\gamma$ from $\partial \Sigma^-$ to $\partial \Sigma^+$. We want to place pairwise disjoint balls on $\gamma$. We have that $d((u \circ \gamma)(1), A^-) \geq \delta$, so by the intermediate value theorem for each $0 \leq n < N$ there is a $t_n \in [0, 1]$ s.t. $d((u \circ \gamma)(t_n), A^-) = (2n + 1)r_0$. Also there is a $t_N$ such that $d((u \circ \gamma)(t_N), A^-) = 2Nr_0 + \delta_r/2$. Picking $x_n = (u \circ \gamma)(t_n)$ for $0 \leq n \leq N$ and using the triangle inequality we have that $\{B_{r_0}(x_n) \mid 0 \leq n < N\} \cup \{B_{\delta_r/2}(x_N)\}$ are pairwise disjoint. See Figure 13. Applying Lemma 5.7 to each $B_r(x_n)$ we get the desired result. □

**Proposition 5.9.** *Take $X = B_{dpq}$ and let $L = T(x)$ be an almost toric fibre (with respect to the almost toric structure described in Section 2.3) with $x \notin R$, where $R$ is defined in (8). Then we have*

$$e(T(x)) = x_1.$$

---

[17]Typically, $r_0$ will be the lower bound on the injectivity radius in Definition 5.2. If $J$ is chosen $\omega$-compatible, we can choose $C = \pi/4$.

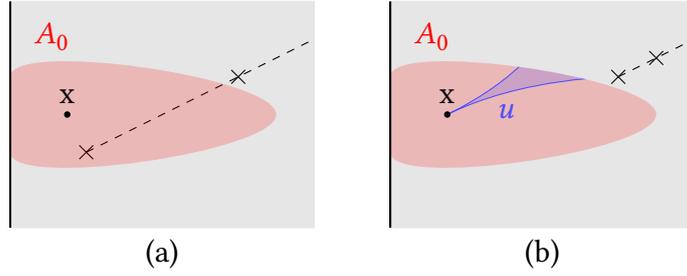

Figure 14: (a) Area $A_0$ obtained from Step 1 is obstructed by node and branch cut. (b) $A_0$ is cleared by a nodal slide. Now $\pi^{-1}(A_0) \cong \mu_0^{-1}(A_0)$. Any J-curve $u$ escaping $\pi^{-1}(A_0)$ must have sufficiently large area.

*Proof. Step 1:* Recall that, as a set, we have $\Delta_{dpq} = \mathbb{R}_{\geq 0} \times \mathbb{R}$. Without any nodes, this is the toric moment image of $\mathbb{C} \times T^*S^1$ under the moment map

$$\mu_0 : \mathbb{C} \times T^*S^1 \to \mathbb{R}_{\geq 0} \times \mathbb{R}$$
$$(z, p, \theta) \mapsto (\pi|z|^2, p),$$

where $\mathbb{C} \times T^*S^1$ is equipped with the symplectic structure $\omega_0 = \omega_{\mathbb{C}} \oplus \omega_{T^*S^1}$, the two direct summands being the standard symplectic structures on $\mathbb{C}$ and $T^*S^1$ respectively. We also equip $\mathbb{C} \times T^*S^1$ with the complex structure $J_0 = i \oplus i$. The second factor of $J_0$ comes from the identification $T^*S^1 = \mathbb{C}/\mathbb{Z}$, where $\mathbb{Z}$ acts on $\mathbb{C}$ by translation along the real axis. This complex structure is compatible with $\omega_0$ and geometrically bounded with respect to the induced metric, which is the Riemannian metric given by

$$d((z_1, p_1, \theta_1), (z_2, p_2, \theta_2))^2 = |z_2 - z_1|^2 + (p_2 - p_1)^2 + \min_{k \in \mathbb{Z}}\{(\theta_2 - \theta_1 + 2\pi k)^2\}$$

for $(z_i, p_i, \theta_i) \in \mathbb{C} \times T^*S^1$.

Let $u, v \in \operatorname{im} \mu_0$. Then we have

$$d(\mu_0^{-1}(u), \mu_0^{-1}(v))^2 \geq (\sqrt{v_1} - \sqrt{u_1})^2 + (v_2 - u_2)^2 =: d_{\mu_0}(u, v)^2, \qquad (15)$$

where $d_{\mu_0}$ defines a metric on $\operatorname{im} \mu_0$.

Let $C, r_0$ be the constants coming from Lemma 5.7 for $(\mathbb{C} \times T^*S^1, \omega_0, J_0)$. Take $\delta > 0$ such that the expression $C(Nr_0^2 + \delta_r^2/4)$ in Corollary 5.8 is greater than $x_1$. Set $A_0^- = \{x\}$,

$$A_0 = \{y \in \operatorname{im} \mu_0 \mid d_{\mu_0}(x, y) \leq \delta\}, \qquad (16)$$
$$A_0^+ = \{y \in \operatorname{im} \mu_0 \mid d_{\mu_0}(x, y) = \delta\}. \qquad (17)$$

*Step 2:* Let us now return to $B_{dpq}$. Since $x \notin R$ and $A_0 \subset \mathbb{R}_{\geq 0} \times \mathbb{R}$ is compact, we can perform a nodal slide as in Figure 14 whose support does not contain $x$ to remove all nodes and branch cuts from $A_0$, see also Lemma 2.4. Let $\pi\colon B_{dpq} \to \Delta_{dpq}$ be the almost toric moment map associated to this almost toric fibration. Since $A_0$ does not contain any nodes or branch cuts, the restriction $\pi|_{\pi^{-1}(A_0)}\colon \pi^{-1}(A_0) \to A_0$ is a toric moment map and thus, using Delzant's theorem [14, Theorem 2.1], we identify $\pi^{-1}(A_0) \subset B_{dpq}$ with $\mu_0^{-1}(A_0) \subset \mathbb{C}\times T^*S^1$ by an equivariant symplectomorphism. Let $J_{B_{dpq}}$ be a geometrically bounded almost complex structure on $B_{dpq}$, which exists by Remark 5.3. Let $K$ be a compact set whose interior contains $\pi^{-1}(A_0)$. By Remark 5.4 there exists a geometrically bounded almost complex structure $J$ on $B_{dpq}$ such that $J = J_{B_{dpq}}$ on $B_{dpq} \setminus K$ and $J = J_0$ on $\mu^{-1}(A_0)$. Set $A^\pm = \pi^{-1}(A_0^\pm)$, $A = \pi^{-1}(A_0)$ and suppose $u\colon (D, \partial D) \to (B_{dpq}, T(x))$ is a $J$-holomorphic disk.

If $u(D) \subset A$, then, under the identification $A = \pi^{-1}(A_0) \cong \mu_0^{-1}(A_0) \subset \mathbb{C}\times T^*S^1$, it defines a relative homology class $[u] \in H_2(\mathbb{C}\times T^*S^1, T(x))$. Therefore $[u]$ is a positive multiple of the generator $D_0$ of $H_2(\mathbb{C}\times T^*S^1, T(x))$, which satisfies $\int_{D_0} \omega_0 = x_1$.

If $u(D)$ is not contained in $A$, then $u(D) \cap A^+$ is non-empty. Let $\Sigma \subset D$ be the connected component of $u^{-1}(A)$ containing $\partial D$. Using Corollary 5.8 and (15) for $A^\pm$ and $u|_\Sigma$, we get

$$\int_u \omega \geq \int_{u|_\Sigma} \omega \geq C(Nr_0^2 + \delta_r^2/4) \geq x_1.$$

In both cases, we conclude that $\int_u \omega \geq x_1$, so $\sigma_D(B_{dpq}, T(x), J) \geq x_1$.

Since $H_2(B_{dpq})$ is generated by Lagrangian spheres (by [21, Lemma 7.11]), there are no non-constant $J$-holomorphic spheres, and $\sigma_S(B_{dpq}, J) = \infty$.

Therefore Theorem 5.6 and Lemma 4.1 allows us to conclude. □

As for the upper bound in Proposition 4.3, there is an analogous statement for small enough tori in partial $B_{dpq}$-embeddings (Definition 2.9).

**Proposition 5.10.** *Let $\varphi\colon B_{dpq} \supset U \to (X, \omega)$ be a partial $B_{dpq}$-embedding into a geometrically bounded symplectic manifold $(X, \omega)$. Then there exists a neighbourhood $V \subset \Delta_{dpq}$ of the origin such that for all $x \in V \setminus R$ we have $e(\varphi(T(x))) = x_1$.*

*Proof.* We proceed in a similar way as in Proposition 5.9: Set $A_0(x) = \{y \in \operatorname{im} \mu_0 \mid d_{\mu_0}(x, y) \leq \delta\}$ as defined in the proof of Proposition 5.9 and take $V \subset \Delta_{dpq}$ a neighbourhood of the origin such that $\pi^{-1}(A_0(x))$ is contained in $U$ for all $x \in V$. Assume additionally that $\varphi(\pi^{-1}(\mathscr{A}_0))$, where $\mathscr{A}_0 = \bigcup_{x \in V} A_0(x)$, is contained in the interior of some compact set $K \subset X$. This can be achieved by further restricting $V$.

Let $J_X$ be a geometrically bounded almost complex structure on $X$, and $J_0$ the almost complex structure on $\pi^{-1}(\mathscr{A}_0)$ obtained by identifying $\pi^{-1}(\mathscr{A}_0)$ with $\mu_0^{-1}(\mathscr{A}_0)$. Let $J$ be a geometrically bounded almost complex structure on $X$ with $J = \varphi_* J_0$ on $\varphi(\pi^{-1}(\mathscr{A}_0))$ and $J = J_X$ on $X \setminus K$. For any fixed $x \in V \setminus R$, we conclude as in Proposition 5.9 that any $J$-disk $u\colon (D, \partial D) \to (X, \varphi(T(x)))$ has at least area $x_1$, so $\sigma_D(X, \varphi(T(x)), J) \geq x_1$.

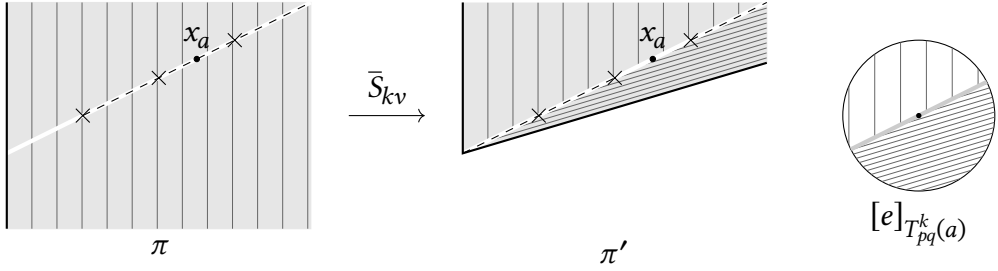

Figure 15: Image of the almost toric moment maps $\pi$ and $\pi'$, before and after changing the branch cut. To the right are the level sets of the displacement energy germ of the fibre $T^k_{pq}(a) = T(x_a)$ as in (19).

Unlike $B_{dpq}$, the space $X$ might contain $J$-spheres. In order to deduce the claim from Theorem 5.6, we additionally require that $x_1 < \sigma_S(X, J)$, which is done by further restricting $V$. By Theorem 5.6 we get $e(\varphi(T(x))) \geq x_1$.

Restricting $V$ once more, we use Proposition 4.3 to get $e(\varphi(T(x))) \leq x_1$. □

# 6 Proof of main theorems

Let $d, p, q \in \mathbb{N}$ be integers such that $d \geq 1$ and $p, q$ coprime. We think of $q$ as an element of $\mathbb{Z}_p/\pm 1$.

Proposition 5.9 yields the first step of the outline given in Section 3.5. To compute the displacement energy germ $[e]_L$ of $L = T^k_{pq}(a)$ (see Definition 2.8) on an open dense subset, we adapt the method in Section 3.5 to the case at hand.

**Theorem 6.1.** *The displacement energy germ of $T^k_{pq}(a) \subset B_{dpq}$ satisfies*

$$[e]_{T^k_{pq}(a)}(b) \sim a + \min\{b_1, b_1(1-kpq) + b_2 kp^2\} \qquad (18)$$

*for some identification $H^1(L; \mathbb{R}) \cong \mathbb{R}^2 = \{(b_1, b_2)\}$.*

*Proof.* Let $\pi \colon B_{dpq} \to \Delta_{dpq}$ be an almost toric moment map with $k$ nodes to the left of the point $x_a = \left(a, \frac{aq}{p}\right)$. Then the almost toric fibre of $\pi$ over $x_a$ is $T^k_{pq}(a)$ by Definition 2.8. Note that if $k \neq 0$, then $x_a$ lies on a branch cut of $\Delta_{dpq}$, as on the left-hand side of Figure 15. To apply Proposition 3.15 we perform changes in branch cut to the $k$ nodes to the left of $x_a$. As a result, we obtain a new almost toric moment map $\pi' \colon B_{dpq} \to \Delta'$, see the right-hand side of Figure 15. By Proposition 2.2, its image $\Delta'$ is given by $\bar{S}^k_v(\Delta_{dpq}) = \bar{S}_{kv}(\Delta_{dpq})$ with $v = (-p, -q)$. In matrix form, the corresponding shear map is given by

$$S_{kv} = \begin{pmatrix} kpq + 1 & -kp^2 \\ kq^2 & 1 - kpq \end{pmatrix}, \quad v = (-p, -q). \qquad (19)$$

Furthermore, the respective almost toric fibres are identified as follows,

$$(\pi')^{-1}(\bar{S}_{kv}(x)) = \pi^{-1}(x), \quad x \in \text{int } \Delta. \tag{20}$$

By Proposition 3.15, we obtain a versal deformation of $L = T_{pq}^k(a)$ by setting

$$v_L(b) = (\pi')^{-1}(x_a + b), \quad b \in V \subset H^1(L; \mathbb{R}) \cong \mathbb{R}^2, \tag{21}$$

where $V$ is a small enough neighbourhood of the origin. As described in Section 3.4, we chose the identification $H^1(L; \mathbb{R}) \cong \mathbb{R}^2$ induced by the identification $H_1(L) \cong \mathbb{Z}^2$ via the closed orbits of the local $T^2$-moment map defined by $\pi'$ in a neighbourhood of $x_a$. Using (20), we can apply Proposition 5.9 to determine the displacement energy of all fibres $(\pi')^{-1}(x)$ for all $x \in \Delta' \setminus R$. In light of the fact that (21) is a versal deformation, this proves the claim. □

Now let $\varphi: B_{dpq} \supset U \to (X, \omega)$ be a partial $B_{dpq}$-embedding in the sense of Definition 2.9. In Proposition 5.10, we have proved that $e(\varphi(T(x))) = e(T(x)) = x_1$ for small enough $x \in \Delta_{dpq}$. This means that the same method as in the proof of Theorem 6.1 can be used to determine the displacement energy germ of such $\varphi(T(x)) \subset X$.

**Theorem 6.2.** *Let $\varphi: B_{dpq} \supset U \to (X, \omega)$ be a partial $B_{dpq}$-embedding into a geometrically bounded symplectic manifold $(X, \omega)$. Then there exists $a_0 > 0$ such that for all $a < a_0$, the displacement energy germ of $\varphi(T_{pq}^k(a)) \subset X$ satisfies*

$$[e]_{\varphi(T_{pq}^k(a))}(b) \sim a + \min\{b_1, b_1(1 - kpq) + b_2 kp^2\}. \tag{22}$$

*Proof.* Let $V \subset \Delta_{dpq}$ be the subset given by Proposition 5.10 such that $e(\varphi(T(x))) = x_1$ for $x \in V$. Let $a_0$ be the supremum over all $a$ such that $x_a = \left(a, \frac{aq}{p}\right) \in V \setminus R$. Picking $a < a_0$ and setting $L = T_{pq}^k(a)$ we construct the versal deformation $v_L$ in the same manner as in Theorem 6.1. Now $v_{\varphi(L)} = \varphi \circ v_L \circ \varphi^*$ is a versal deformation of $\varphi(L)$, where $\varphi^*$ is the induced map on cohomology. As $\varphi$ preserves the displacement energy of $T(x)$ for $x \in V \setminus R$ by Propositions 5.9 and 5.10, we have $[e]_L \sim [e]_{\varphi(L)} \circ \varphi^*$. □

Using Theorems 6.1 and 6.2, we can prove Theorems A and B:

*Proof of Theorems A and B.* We show that if $[e]_{T_{pq}^k(a)} \circ \Phi \sim [e]_{T_{p'q'}^{k'}(a')}$ for some $\Phi \in \text{GL}(2; \mathbb{Z})$, then $a = a'$ and either $k = k' = 0$ or $(k, p, q) = (k', p', q')$. The value of $a$ can be directly read of off (18) and (22), so $a = a'$. Equations (18) and (22) are both given as the minimum of two linear functions $\alpha_{pq}^k, \beta_{pq}^k$, where

$$\alpha_{pq}^{ka}(b) = a + \det\left(b, \begin{pmatrix} 0 \\ 1 \end{pmatrix}\right) \quad \beta_p^k q(b) = a - \det\left(b, \begin{pmatrix} kp^2 \\ kpq - 1 \end{pmatrix}\right).$$

The set $\{\alpha_{pq}^k = \beta_{pq}^k\}$ is given by the space spanned by the vector $v = (p, q)$. Denote by $u_1 = \begin{pmatrix} 0 \\ 1 \end{pmatrix}$ and $u_2 = \begin{pmatrix} kp^2 \\ kpq - 1 \end{pmatrix}$ the vectors pointing along the level sets of $\alpha_{pq}^k$ and $\beta_{pq}^k$ respectively, and define the vectors $v', u_1', u_2'$ similarly.

Similarly to Remark 2.14, the integral affine angle between $v$ and $u_1$ (or $u_2$) determines $(p, q) \in \mathbb{N}_{>0} \times \mathbb{Z}_p/\pm 1$, and $u_1 + u_2 = k(p, q)$, determining $k$. So if $\min\{\alpha_{pq}^k, \beta_{pq}^k\} \circ \Phi \sim \min\{\alpha_{p'q'}^{k'}, \beta_{p'q'}^{k'}\}$, then $\Phi(v') = v$, $\Phi(\{u_1', u_2'\}) = \{u_1, u_2\}$ and by the above $(k, p, q) = (k', p', q')$.

We show that $T_{pq}^k(a) \sim T_{pq}^{k'}(a)$. Denote by $\pi_k, \pi_{k'}$ the almost toric moment maps used in Definition 2.8 for $T_{pq}^k(a), T_{pq}^{k'}(a)$ respectively. We have $T_{pq}^k(a) = \pi_k^{-1}(x_a)$ and $T_{pq}^{k'}(a) = \pi_{k'}^{-1}(x_a)$. Suppose that $\varepsilon > 0$ and take the Lagrangian isotopies $\Lambda_t^k = \pi_k^{-1}(x_a + (0, t\varepsilon))$ and $\Lambda_t^{k'} = \pi_{k'}^{-1}(x_a + (0, t\varepsilon))$. By Proposition 3.15, $\text{Flux}(\Lambda^k) = \text{Flux}(\Lambda^{k'}) = (0, \varepsilon)$. Using a nodal slide as in Lemma 2.4 we may identify $\Lambda_1^k$ and $\Lambda_1^{k'}$. Concatenating $\Lambda^k$ and the inverse of $\Lambda^{k'}$ we get a Lagrangian isotopy $\Lambda$ from $T_{pq}^k(a)$ to $T_{pq}^{k'}(a)$ with $\text{Flux}(\Lambda) = \text{Flux}(\Lambda^k) - \text{Flux}(\Lambda^{k'}) = 0$. The Maslov class is preserved under Lagrangian isotopies, and $\text{Flux}(\Lambda) = 0$ means that the area class is preserved as well. $\square$

*Remark* 6.3. As mentioned in Remark 1.1, some special cases of Theorem A can be deduced from [25]. Indeed,

1. for $p = q = 1$ and arbitrary $d$, the space $B_{d,1,1}$ appears as Milnor fibre of the smoothing of the $A_{d-1}$-singularity (see also [21, Section 7.3]) and [25, Lemma 2.19] yields a count of Maslov index two $J$-holomorphic disks with boundary on the tori $T_{1,1}^k$. This count depends on $k$ and, by monotonicity of the tori, is a symplectic invariant, see [17]. Thus [25, Lemma 2.19] distinguishes the tori.

2. for $d = 1, p \neq 1$ and $(p, q) \neq (1, 0)$, the torus $T_{pq}^1(a)$ is non-displaceable by [25, Proposition 3.6]. Since $T_{pq}^0$ is displaceable by the methods used in Section 4, we deduce $T_{pq}^0 \not\cong T_{pq}^1$.

The special case of $d = 2, p = 1, q = 0$ corresponds to $B_{2,1,0} \cong T^*S^2$ and thus follows from earlier work [2], see also the discussion in Example 1.2. Our method of proof is different from the one used in [25]. Instead of a count of Maslov index two $J$-holomorphic disks, we rely on the displacement energy germ to distinguish tori. The upshot of the method used in [25] is that it yields non-displaceability results via computing the Floer homology. Our methods, on the other hand can be applied to non-monotone tori, which is crucial in Theorem B. Because of bubbling, the count of Maslov index two $J$-holomorphic disks depends on the choice of $J$ and does not yield an invariant in that case.

# 7 Examples

Along with Corollaries 1.3, 1.4 and 1.6 given in the introduction, we have the following:

As discussed in [21, Remark 7.10], if $d > 1$, the space $B_{dpq}$ contains a Lagrangian $S^2$. Together with Corollary 1.3, this means every time we have a partial $B_{dpq}$-embedding with $d > 1$, we also get a partial $B_{2,1,0}$-embedding "for free". In particular a Polterovich torus $T^2_{1,0} \subset B_{2,1,0}$ is contained in any partial $B_{dpq}$-embedding with $d > 1$. By Theorem B small enough such tori are distinct from the $T^k_{pq}$ coming from the same embedding.

## 7.1 Stretching an edge

**Proposition 7.1.** *Let $\Delta \subset \mathbb{R}^2$ be a compact (almost) toric base diagram with finitely many nodes, and $(X, \omega)$ the corresponding (almost) toric manifold. Then there is a sequence $(d_i, p_i, q_i)$ with $p_i \to \infty$ such that there is a partial $B_{d_i, p_i, q_i}$-embedding into M.*

We can assume, after performing nodal trades, that every vertex of $\Delta$ has one or multiple corresponding nodes and, by nodal slides, that all the branch cuts are negligibly small.

Recall from Definition 2.11 that an ATF-corner $u$ consists of a pair $u = (x, dv)$ given by a vertex $x$ and a multiple $dv$ of the branch cut direction.

**Definition 7.2.** The **mutation at** $(x, dv)$ is given by performing a change of branch cut at all nodes lying on the branch cut $v$ as in Proposition 2.2, followed by nodal slides to make the newly created branch cut negligibly small.

Also, recall that the piece-wise linear map associated to the change of branch cut is given by $\bar{S}_{(x,dv)}$, as defined in (7).

The idea of the proof of Proposition 7.1 is to perform an infinite sequence of mutations on $\Delta$, to obtain a sequence of partial $B_{d_i, p_i, q_i}$-embeddings. In order to avoid some cyclical behaviour in the sequence $(d_i, p_i, q_i)$, we fix an edge and perform iterated mutations at its left end. Under this operation the integral affine length of the edge in question increases.

*Remark* 7.3. Let us briefly describe what the algorithm in the proof of Proposition 7.1 yields for the case of $X = \mathbb{C}P^2$. Vianna [34] proved that the set of ATF-base diagrams of $\mathbb{C}P^2$ obtained from the toric Delzant polytope by mutations is, up to nodal slides, in bijection with so-called Markov triples, i.e. solutions to the Markov equation. In the case of $\mathbb{C}P^2$ it follows immediately from the definition of the sequence of mutations used in the proof of Proposition 7.1 that, starting from the standard toric structure, this sequence corresponds to the unique path of Markov triples in the Markov tree that contain a one. If we start at another Markov triple, it corresponds to the path in the Markov tree obtained by fixing one number in the Markov triple.

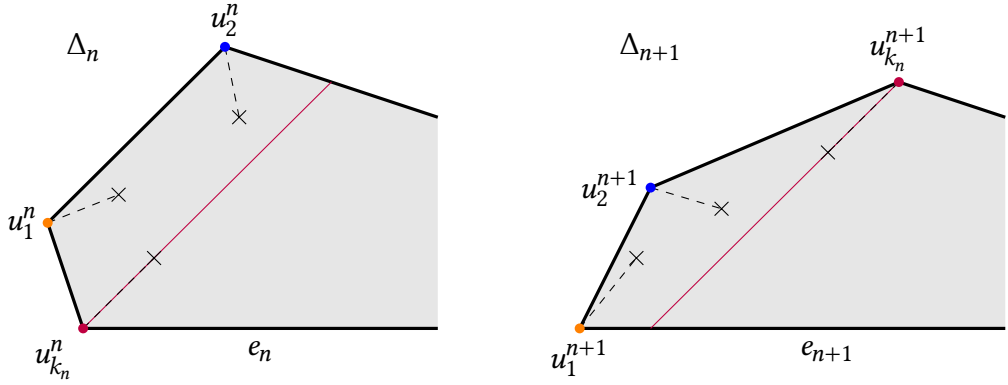

Figure 16: Moving from $\Delta_n$ to $\Delta_{n+1}$ by performing a mutation at $u^n_{k_n}$.

For more details on this, and in particular the combinatorics of the sequences involved, we refer to [21, Section 8.4, Appendix I].

*Proof of Proposition 7.1.* To simplify, we slightly abuse notation and denote by $u$ the ATF-corner as well as the corresponding vertex $x$. Throughout the proof we assume, possibly after performing nodal slides, that all branch cuts are negligibly small.

*Step 1*: Set $\Delta_0 = \Delta$ and pick an edge $e_0$ of $\Delta_0$. We may assume $e_0 \subset \mathbb{R} \times \{0\}$ and $\Delta_0 \subset \mathbb{R} \times \mathbb{R}_{\geq 0}$. We will elongate $e_0$ by repeated iterations at an adjacent ATF-corner. In other words, we will define a sequence $\{\Delta_n\}_{n \in \mathbb{N}}$ of almost toric base diagrams obtained by mutation from $\Delta_0$, which contain a distinguished edge $e_n \subset \Delta_n$.

Let $u^0_1, ..., u^0_{n_c}$ be the ATF-corners of $\Delta_0$, where $n_c$ is the number of ATF-corners. Let $\Delta_{n+1} = \bar{S}_{u^n_{k_n}} \Delta_n$ be the result of mutating $\Delta_n$ at the ATF-corner $u^n_{k_n}$ next to $e_n$ in the clockwise direction and if $e_n$ has multiple non-isolated ATF-corners at this vertex, we pick $u^n_{k_n} = (x, v)$ for which $v$ forms the smallest angle with $e_n$. Let $u^{n+1}_1, ..., u^{n+1}_{n_c}$ be the corners of $\Delta_{n+1}$, such that for $i \neq k_n$ we have $\bar{S}_{u^n_{k_n}} u^n_i = u^{n+1}_i$ and such that for $i = k_n$, $u^{n+1}_{k_n}$ is the new ATF-corner created by the mutation at $u^n_{k_n}$.[18] In other words, the elements in the set $\{1, ..., n_c\}$ can be thought of as labels of ATF-corners, where every ATF-corner in $\Delta_{n+1} = \bar{S}_{u^n_{k_n}} \Delta_n$ inherits the label of its preimage in $\Delta_n$, except for the newly created ATF-corner in $\Delta_{n+1}$, which inherits the label of the ATF-corner in $\Delta_n$ at which mutation was performed. See Figure 16, where we have indicated the labels by colours. Note that $e_n$ lies in the half plane fixed by $\bar{S}_{u^n_{k_n}}$. Take $e_{n+1}$ to be the edge of $\Delta_{n+1}$ such that $e_{n+1} \supset e_n$.

*Step 2*: If $u^n_{k_n}$ is not isolated, we have $\ell(e_{n+1}) = \ell(e_n)$, where $\ell(e)$ denotes the affine length of an edge $e$, otherwise $\ell(e_{n+1}) > \ell(e_n)$. Since the first case can only appear finitely many times successively, the monotone sequence $\ell(e_n)$ has a strictly monotone subsequence. It is also bounded from above by the affine length of $\partial \Delta_0$. Define $\ell_n = \ell(e_n) - \ell(e_0)$ and denote by $\ell_\infty$ the limit of $\ell_n$ as $n \to \infty$.

---

[18]After a mutation we might loose an ATF-corner if the eigenline of $u^n_{k_n}$ intersects $\partial \Delta_n$ in another ATF-corner with the same eigenline. We can assume without loss of generality that this does not happen. Indeed, the sequence $\Delta_n$ must attain its minimum number of ATF-corners, so it will be constant for all $n \geq N$ for some $N \in \mathbb{Z}$ and we then take $\Delta_0 = \Delta_N$.

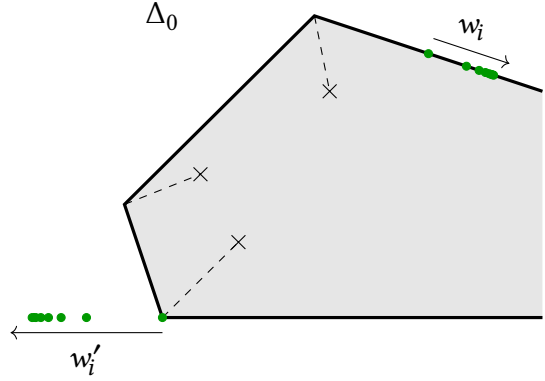

Figure 17: Artist's impression of the convergence of $w_i$ and $w_i'$.

*Step 3*: The sequence $k_n \in \{1, ..., n_c\}$ of labels contains some elements of $\{1, ..., n_c\}$ infinitely many times. We may assume without loss of generality that all the elements which appear, appear infinitely many times. If this is not the case, we can choose $N$ such that this holds for $k_n$ with $n \geq N$ and then set $\Delta_0 = \Delta_N$. Furthermore, we can assume that the labels appearing infinitely many times are $\{1, ..., n_c' \leq n_c\}$, which can be achieved by relabelling.

Note that the vertex $u_{k_0}^0 \in e_0$ at which we have performed the first mutation is fixed by all subsequent $\bar{S}_{u_{k_n}^n}$, so it makes sense to speak of $u_{k_0}^0 \in \partial \Delta_n$ for all $n$. Denote by $g_n(x)$ the integral affine length of the portion of $\partial \Delta_n$ contained between $u_{k_0}^0$ and $x \in \partial \Delta_n$ in the clockwise direction. Note that we have $g_n(x) = g_{n+1}(\bar{S}_{u_{k_n}^n} x)$, since integral linear transformations preserve the affine length.

Define $a_n = \max_{k \in \{1,...,n_c'\}} \{g_n(u_k^n)\}$, which is the distance from $u_{k_0}^0$ to the ATF-corner furthest away in terms of $g_n(\cdot)$ appearing in infinitely many mutations. We have $a_n \leq a_{n+1}$, since $g_n(u_k^n) = g_{n+1}(u_k^{n+1})$ if $k \neq k_n$ and $g_n(u_{k_n}^n) < g_{n+1}(u_{k_n}^{n+1})$ otherwise. Also $a_n \geq \ell_n = g_n(u_{k_n}^n)$ and $a_n \leq \ell_m$ for some $m > n$. The latter claim follows from the fact that for $k \in \{1, ..., n_c'\}$ realizing the maximum in the definition of $a_n$ we must have $k = k_m$ for some $m > n$, as $k$ appears infinitely many times in $\{k_m\}_{m \in \mathbb{N}}$, and thus $a_n = g_n(u_{k_m}^n) \leq g_m(u_{k_m}^m) = \ell_m$. In particular, $a_n \to \ell_\infty$ and we can choose a strictly monotone subsequence $a_{n_i}$ with $a_{n_i} > a_{n_i - 1}$.

*Step 4*: Now chose the natural number $k$ realizing the maximum in the definition of $a_{n_i}$ and define $w_i = u_k^{n_i}$. Actually it follows from the definitions that $k = k_{n_i - 1}$, since $a_{n_i} > a_{n_i - 1}$ means that the maximum of $a_{n_i}$ must be created by the newly created ATF-corner in $\Delta_{n_i}$, which is $u_{k_{n_i-1}}^{n_i}$. We have $w_i \in \partial \Delta_0$: Indeed, the maps $\bar{S}_{u_{k_n}^n} : \Delta_n \to \Delta_{n+1}$ divide $\partial \Delta_n$ into two parts: If $x \in \partial \Delta_n$ satisfies $g_n(x) \in \left[g_n(u_{k_n}^n), g_{n+1}(u_{k_n}^{n+1})\right]$, i.e. $x$

lies in the interval in $\partial \Delta_n$ from $u^n_{k_n}$ to $u^{n+1}_{k_n}$ in the clockwise direction, we have $\bar{S}_{u^n_{k_n}} x = S_{u^n_{k_n}} x$. Otherwise, $x$ is in the half plane fixed by $\bar{S}_{u^n_{k_n}}$. Since $g_{n_i}(w_i) > g_n(u^n_k)$ for all $k \in \{1, ..., n_c\}$ and $n < n_i$, the point $w_i$ lies in the intersection of the half planes that are fixed by all the mutations that are applied in the sequence from $\Delta_0$ to $\Delta_{n_i}$.

Define $w'_i = u^{n_i-1}_{k_{n_i-1}}$, which is the ATF-corner at which we mutated to obtain $\Delta_{n_i}$. Denote by $w_\infty \in \partial \Delta_0$ the limit point of $\{w_i\}_{i \in \mathbb{N}}$ and denote by $w'_\infty$ the limit point of $\{w'_i\}_{i \in \mathbb{N}}$, contained in the line spanned by $e_0$. See Figure 17. Define $v_i$ to be the primitive vector in the ATF-corner $w_i = (x, dv_i)$ pointing along $w'_i - w_i$, and $v_\infty = w'_\infty - w_\infty$. Denote by $f$ the edge of $\Delta_0$ containing infinitely many $w_i$, as well as $w_\infty$, and define $x_f$ to be the primitive vector pointing along $f$ in the counter-clockwise direction. Analogously, define $x_e$ to be the primitive vector pointing along $e_0$ in the counter-clockwise direction. Then for $i$ large enough $v_{i+1} = \eta v_i + \lambda x_f - \mu x_e$ with $\eta, \lambda, \mu > 0$. We have that $\det(x_f, v_i) > 0$, by convexity of the ATF-base diagrams $\Delta_n$, and also $\det(x_e, v_i) < 0$. So

$$\det(v_{i+1}, v_i) = \lambda \det(x_f, v_i) - \mu \det(x_e, v_i) > 0 .$$

Denote the projection onto the unit circle by $\pi \colon \mathbb{R}^2 \setminus \{0\} \to S^1$. Since $\det(x_e, v_i) < 0$, the points $\pi(v_i)$ are contained in a half-circle, and thus they converge strictly monotonically to $\pi(v_\infty)$. Since $v_i$ is primitive, the Euclidean length $\|v_i\|$ satisfies $\|v_i\| \to \infty$.

The vector $v_\infty$ is not parallel to $x_f$, since it would contradict the convergence of $w_i$ if it were parallel. So $p_i = \det(x_f, v_i) \to \infty$ and $w_i$ is an ATF-corner of type $(d_i, p_i, q_i)$ in $\Delta_{n_i}$.

Using Lemma 2.13 we get a partial $B_{d_i p_i q_i}$-embedding for every $i \in \mathbb{Z}$. □

*Remark* 7.4. Based on simulations,[19] the situation described in the proof is always very simple: We have $n'_c = 2$ so that the number of ATF-corners appearing infinitely many times is two, so we always swap back and forth, as in sequences of mutations used in [18]. This would also allow to prevent the sizes of embeddings in Theorem C to become arbitrarily small.

---

[19]Done using https://github.com/schmijoe/delzant.